\documentclass[11pt]{elsarticle}

\usepackage{lineno,hyperref}
\usepackage[margin=3.9cm]{geometry} 

\journal{}
\usepackage[font=small]{caption}
\usepackage{amsmath}
\usepackage{amssymb}

\usepackage{hyperref}

\usepackage{tikz}
\usepackage{pgfplots,xcolor,soul}

\usepackage{algorithm}
\usepackage{algorithmicx}
\usepackage[noend]{algpseudocode}

\usepackage{stmaryrd}

\usepackage{longtable}

\usepackage{booktabs}

\usepackage{multirow}

\usepackage{subcaption}

\usepackage[title]{appendix}

\newtheorem{definition}{Definition}
\newtheorem{theorem}{Theorem}
\newtheorem{lemma}{Lemma}

\newenvironment{proof}[1][Proof]{\textbf{#1:}\\}{\hfill\rule{1.5mm}{1.5mm}\\}

\newcommand{\init}{init}
\newcommand{\new}{new}
\newcommand{\rec}{rec}

\newcommand{\mitenchange}[1]{{#1}}

\allowdisplaybreaks

\begin{document}

\begin{frontmatter}

\title{Exact Lexicographic Scheduling and Approximate Rescheduling}

\author{Dimitrios Letsios$^{a,*}$}
{\fnref{label0}
\corref{cor1}
}

\author{Miten Mistry$^{b,**}$}
{\fnref{label1}
\corref{cor2}
}

\author{Ruth Misener$^{b,**}$}
{\fnref{label1}
\corref{cor2}}

\address[label0]{Department of Informatics; King's College London; United Kingdom}
\address[label1]{Department of Computing; Imperial College London; South Kensington SW7 2AZ; UK}

\cortext[cor1]{\texttt{dimitrios.letsios@kcl.ac.uk}}
\cortext[cor2]{\texttt{\{miten.mistry11, r.misener\}@imperial.ac.uk}; Tel: +44 (0) 20759 48315}

\begin{abstract}
In industrial resource allocation problems, an initial planning stage may solve a nominal problem instance and a subsequent recovery stage may intervene to repair inefficiencies and infeasibilities due to uncertainty, e.g.\ machine failures and job processing time variations.  
In this context, we investigate the minimum makespan scheduling problem, a.k.a.\ $P||C_{\max}$, under uncertainty. 
We propose a two-stage robust scheduling approach where first-stage decisions are computed with exact lexicographic scheduling and second-stage decisions are derived using approximate rescheduling.
We explore recovery strategies accounting for planning decisions and constrained by limited permitted deviations from the original schedule.
Our approach is substantiated analytically, with a price of robustness characterization parameterized by the degree of uncertainty, and numerically.
This analysis 
is based on optimal substructure imposed by lexicographic optimality. 
Thus, lexicographic optimization enables more efficient rescheduling. 
Further, we revisit state-of-the-art exact lexicographic optimization methods and propose a 
lexicographic 
branch-and-bound algorithm whose performance is validated computationally.
\end{abstract}


\end{frontmatter}



\section{Introduction}
\label{Section:Introduction}

\begin{figure}
\begin{center}
\includegraphics{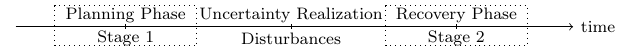}
\caption{Recoverable robustness setting}
\label{Figure:Recoverable_Robustness_Model}
\end{center}
\end{figure}

Motivated by industrial resource allocation problems, we consider scheduling under uncertainty, e.g.\
a machine may unexpectedly fail, a client may suddenly cancel a job, or 
jobs 
are completed earlier than expected.
We focus on robust scheduling, i.e.\ hedge against worst-case realizations of imprecise parameter values, such as job processing times and number of available machines, lying in well-defined 
uncertainty sets \citep{BenTal2009,Bertsimas2011,Goerigk2016, Wiebe2018}.
Because static robust optimization 
may produce conservative solutions compared to ones obtained with perfect knowledge \citep{Soyster1973}, 
we investigate two-stage robust optimization with recovery \citep{BenTal2004,Bertsimas2010,Bertsimas2018,Hanasusanto2015,Liebchen2009}. 
As shown in Figure \ref{Figure:Recoverable_Robustness_Model}, (i) an initial planning stage computes a solution with nominal parameter values and (ii) a subsequent recovery stage modifies the solution once the uncertainty is realized, i.e.\ after the final parameter values become known.

\begin{figure}
\begin{center}
\includegraphics{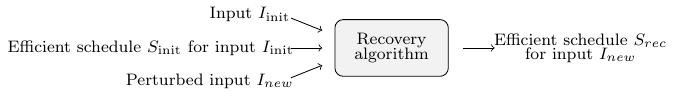}
\caption{Makespan recovery problem}
\label{Figure:Recovery_Makespan_Problem}
\end{center}
\end{figure}

We elaborate on the fundamental makespan scheduling problem, a.k.a.\ $P||C_{\max}$ \citep{Brucker2007,Graham1969,Leung2004}.
With perfect knowledge, an instance $I$ of the problem posits a set $\mathcal{J}$ of jobs, each one associated with processing time $p_j$, a set $\mathcal{M}$ of parallel identical machines and the objective is to construct a non-preemptive schedule $S$ 
of minimum makespan $C_{\max}=\max_{i\in\mathcal{M}}\{C_i\}$, i.e.\ maximum machine completion time. 
In a two-stage setting under uncertainty, the planning stage solves a nominal instance $I_{\init}$ producing a solution $S_{\init}$ and the recovery stage transforms schedule $S_{\init}$ to a new schedule $S_{\new}$ for $I_{\new}$ by repairing inefficiencies and infeasibilities, e.g.\ due to job processing time variations and machine failures, as illustrated in Figure~\ref{Figure:Recovery_Makespan_Problem}.
In the extreme case, 
the recovery stage
may solve $I_{\new}$ from scratch without 
accounting for the first-stage 
decisions in $S_{\init}$.
This level of flexibility may sacrifice benefits due to planning and can be resource-consuming. 
For example, significantly modifying machine schedules may incur substantial communication costs
in distributed computing \citep{yu2007adaptive}.
To mitigate this overhead, we only allow a bounded number of modifications to $S_{\init}$. 
Technically, we distinguish between \emph{binding} and \emph{free} optimization decisions. 
Binding decisions are variable evaluations determined from the initial solution 
after uncertainty realization.
Free decisions are variable evaluations that cannot be determined from the initial solution, 
but are essential to ensure feasibility.
%
%
%
For instance, 
scheduling a job with a modified processing time is a binding decision because 
the planning stage 
already specifies an assignment.
Assigning a new job after uncertainty realization is a free decision because no assignment is given in 
the planning stage. 
Further, we study rescheduling 
with limited binding decision modifications and thereby stay close to the 
initial solution.
By allowing few modifications, first-stage decisions remain critical. 



A two-stage robust optimization method should specify
(i) 
a way of producing the initial solution $S_{\init}$, and (ii) a recovery strategy for restoring $S_{\init}$ and deriving the recovered solution $S_{\rec}$, 
after uncertainty realization. 
Analyzing a two-stage robust optimization method requires defining (i) the uncertainty set of the problem, and (ii) the investigated performance guarantee.


\paragraph{Uncertainty Set}
The uncertainty set of a robust optimization problem specifies a range of possible values for the uncertain parameters \citep{Kouvelis2013}.
We consider a generalization of well-known $\Gamma$-uncertainty sets, where the final parameter values $\hat{p}_j$
vary in an interval $[p_j^L,p_j^U]$ and at most $k$ parameters can
deviate from their nominal values \citep{Bertsimas2003}.
Here, the uncertainty set is defined by a pair $(k,f)$, where 
$k$ is the maximum number of unstable parameters with respect to perturbation factor $f>1$.
A parameter $p_j>0$ is stable if $p_j/f\leq \hat{p}_j\leq fp_j$ and unstable, otherwise.

\paragraph{Performance Guarantee}
Theoretical performance guarantees are useful for determining when robust optimization methods are efficient \citep{Bertsimas2011,Goerigk2016}. 
Denote by $C(I_{\new})$ the cost, e.g.\ makespan, of a solution obtained by some robust optimization method and by $C^*(I_{\new})$ the cost of an optimal solution obtained with perfect knowledge \citep{Monaci2013}.
We consider the so-called \emph{price of robustness} which is defined as the ratio between the two and
seek a tight, worst-case performance guarantee $\rho=\max_{I_{\new}\in\mathcal I} (C(I_{\new}) / C^*(I_{\new}))$ within the set $\mathcal{I}$ of all problem instances \citep{Bertsimas2004}. 

\paragraph{Related Work}

Prior literature 
shows that scheduling problems become computationally harder after incorporating uncertainty \cite{Gupta2019, Gupta2016, Wiebe2018}.
\citet{Kasperski2014} survey techniques and negative results for robust scheduling, including $P||C_{\max}$ with uncertain job processing times. 
Typically, robustness is achieved by optimizing the worst-case (i) cost or (ii) distance from the achievable optimum with perfect knowledge, over all scenarios in the uncertainty set.
These robust counterparts of standard deterministic optimization problems are often referred to as
\emph{minmax} and \emph{minmax regret}, 
respectively \cite{Kouvelis2013}.

With perfect knowledge, $P||C_{\max}$ is strongly $\mathcal{NP}$-hard, but admits greedy constant-factor approximation algorithms and polynomial-time approximation schemes (PTASs) \cite{Brucker2007, Leung2004}.
When the number of machines is constant, $Pm||C_{\max}$ is weakly $\mathcal{NP}$-hard and has fully polynomial-time approximation schemes.
With budgeted uncertainty, where $B$ bounds the deviation of the sum of processing times from their nominal values, minmax $P||C_{\max}$ admits a PTAS when $B$ is constant and a 3-approximation algorithm for arbitrary $B$ \cite{Bougeret2019}.

To our knowledge, no prior work analyzes the price of robustness for $P||C_{\max}$ under uncertainty.
Determining the price of robustness for fundamental combinatorial optimization problems has been repeatedly posed as an open question by domain experts \cite{Bertsimas2011,Bertsimas2004,Goerigk2016}.
The current manuscript shows that such an analysis provides useful structural properties and quantifies the effect of uncertainty in robust solutions for $P||C_{\max}$. 

\paragraph{Lexicographic Optimization}
\ref{eq:lexopt} is at the core of our two-stage scheduling approach.
\ref{eq:lexopt} is a subclass of multiobjective optimization
minimizing $m$ objective functions $F_1,\ldots,F_m:\mathcal{S}\rightarrow\mathbb{R}_0^+$, in decreasing priority order \citep{Erghott2006,Pardalos2016}. 
In other words, \ref{eq:lexopt} optimizes the highest-rank objective $F_1$, then the second most important objective $F_2$, then the third $F_3$, etc.:
\begin{equation}
\tag{LexOpt}\label{eq:lexopt}
\text{lex}\min\{F_1(S),\ldots,F_m(S):S\in\mathcal{S}\}.
\end{equation}
%
%
There are indications that \ref{eq:lexopt} is useful in optimization under uncertainty.
\ref{eq:lexopt} maintains a good approximate schedule when jobs are added and deleted dynamically \citep{sanders2009online, Skutella2016}.
\ref{eq:lexopt} is also useful for cryptographic systems against attacks \citep{Zufiria2017}.
We consider the \ref{eq:lexopt} scheduling problem $\text{lex}\min\{C_1(S),\ldots,C_m(S):S\in\mathcal{S}\}$ of computing a schedule $S$ with lexicographically minimal machine completion times and show that it enables more efficient two-stage robust scheduling.
That is, we identify robust scheduling as a new \ref{eq:lexopt} application.

Apart from optimization under uncertainty, designing efficient \ref{eq:lexopt} methods is motivated by \ref{eq:lexopt} applications: equitable allocation of a divisible resource \citep{Georgiadis2002,Luss1999}, 
fairness \citep{Bouveret2009},
and exploiting opponent mistakes in game theory
\citep{Nace2007,Schmeidler1969}.
Solution strategies include sequential, weighting, and highest-rank objective methods \citep{Burkard1991,Cramer1979,Erghott2006,Ogryczak1997,Pardalos2016,Sherali1982,Sherali1983}. 
There is work characterizing the convex hull of \ref{eq:lexopt} problems \citep{Adams2016,GUPTE2016150,MULDOON201332}.
Logic-based methods are also applicable \cite{Mistry2018}.


\paragraph{Contributions and Paper Organization}
Our main contribution is a two-stage robust scheduling approach for $P||C_{\max}$ under uncertainty, where first-stage decisions are computed with mixed-integer linear programming and lexicographic optimization, while second-stage decisions are derived using approximation algorithms.
Despite the relevant literature on two-stage robust optimization for various applications \cite{BenTal2004,Chassein2016,CHASSEIN2018423,Liebchen2009}, we are not aware of any work on cornerstone scheduling problems, such as $P||C_{\max}$, combined with a price of robustness characterization.

The manuscript proceeds as follows.
Section 
\ref{Sec:Problem_Definitions} formally defines $P||C_{\max}$, 
\ref{eq:lexopt} scheduling, 
and the considered perturbation types.
Section \ref{Section:Lexicographic_Optimization} develops a branch-and-bound algorithm for the \ref{eq:lexopt} scheduling problem.
Section \ref{Section:Recovery} proposes a recovery strategy and analyzes the performance of the overall two-stage approach theoretically.
Section \ref{Section:Numerical_Results} substantiates the branch-and-bound method with respect to state-of-the-art \ref{eq:lexopt} approaches adapted to \ref{eq:lexopt} scheduling.
Further, Section \ref{Section:Numerical_Results} validates our two-stage method empirically.
Section \ref{Section:Conclusion} concludes.

After proving that the makespan recovery problem is strongly $\mathcal{NP}$-hard, at least as hard as solving the problem with full input knowledge, we elaborate on performance guarantees for two-stage $P||C_{\max}$ under uncertainty.
Technically, we investigate a basic recovery strategy that enforces all available binding decisions and performs only essential actions to regain feasibility.
On the negative side, every recovered solution is a weak approximation if planning produces an arbitrary nominal optimal solution.
Specifically, every recovered solution attains an $\Omega(m)$ price of robustness, even in the case of a single perturbation. 
On the positive side, we obtain significantly better performance guarantees if the initial solution is \ref{eq:lexopt}. 
For a single perturbation, 
planning using \ref{eq:lexopt} ensures a price of robustness equal to 2.
For multiple perturbations, the initial solution can be weakly reoptimizable with a high-degree of uncertainty.
However, we show an asymptotically tight $O(f(1+\frac{k}{m-k})(f+k)(1+\frac{\delta}{m}))$ price of robustness, where $k$ is the number of unstable jobs with respect to a perturbation factor $f$ and $\delta/m$ is the fraction of additional machines, after uncertainty realization. 
This result exploits our uncertainty set structure.
Therefore, when $k$, $f$ and $\delta/m$ are constant, our approach achieves an $O(1)$ price of robustness. 

The main paper includes the proof of Theorem \ref{Theorem:Single_Recovery_Negative} and part of the proofs for Theorems \ref{Theorem:Single_Recovery_Positive}-\ref{Theorem:PG},
which bound the price of robustness of our two-stage robust scheduling approach for $P||C_{max}$ under uncertainty in the case of single and multiple perturbations, respectively, by exploiting the optimal substructure imposed by \ref{eq:lexopt}. 
All remaining proofs are provided in the supplementary material.


\section{Problem Definitions}
\label{Sec:Problem_Definitions}

This section defines the $P||C_{\max}$ problem (Section \ref{Section:Makespan_Problem}), the \ref{eq:lexopt} scheduling problem (Section \ref{Section:Lexicographic_Scheduling_Problem}), and describes the investigated perturbations (Section \ref{Section:Perturbation_Types}).

\subsection{Makespan Scheduling Problem}
\label{Section:Makespan_Problem}

An instance $I$ of the makespan scheduling problem, a.k.a.\ $P||C_{\max}$, is a pair $(m,\mathcal{J})$, where $\mathcal{J}=\{J_1,\ldots,J_n\}$ is a set of $n$ 
jobs, 
with processing times $p_1,\ldots,p_n$, to be executed by a set $\mathcal{M}=\{M_1,\ldots,M_m\}$ of $m$ parallel identical machines.
Job $J_j\in\mathcal{J}$ must be processed by exactly one machine $M_i\in\mathcal{M}$ for $p_j$ units of time non-preemptively, i.e.\ in a single continuous interval without interruptions.
Each machine processes at most one job per time.
The objective is 
to minimize the last machine completion time.
Given a schedule $S$, let $C_{\max}(S)$ and $C_i(S)$ be the makespan and the completion time of machine $M_i\in\mathcal{M}$, respectively, in $S$. 
In the following mixed-integer linear programming (MILP) formulation, 
binary variable $x_{i,j}$ is 1 if job $J_j\in\mathcal{J}$ is executed by machine $M_i\in\mathcal{M}$ and 0, otherwise.
\begin{subequations}
\label{Eq:Makespan_MILP}
\begin{align}
\min_{C_{\max},C_i,x_{i,j}} \quad & C_{\max} \label{Eq:Makespan_Objective} \\ 
& C_{\max} \geq C_i && M_i\in\mathcal{M} \label{Eq:Constraint_Max_Completion} \\
& C_i = \sum_{j=1}^nx_{i,j}\cdot p_j && M_i\in\mathcal{M} \label{Eq:ConstraintMachineCompletionTime} \\
& \sum_{i=1}^m x_{i,j} = 1 && J_j\in\mathcal{J} \label{Eq:ConstraintJobAssignment} \\
& x_{i,j}\in\{0,1\} && J_j\in\mathcal{J}, M_i\in\mathcal{M}.
\end{align}
\end{subequations}
Expression (\ref{Eq:Makespan_Objective}) minimizes makespan.
Constraints (\ref{Eq:Constraint_Max_Completion}) enforce that $C_{\max}=\max_{1\leq i\leq m}\{C_i\}$.
Constraints (\ref{Eq:ConstraintMachineCompletionTime}) ensure that a machine executes at most one job per time.
Constraints (\ref{Eq:ConstraintJobAssignment}) impose that each job is assigned to exactly one machine.

\subsection{LexOpt Scheduling Problem}
\label{Section:Lexicographic_Scheduling_Problem}


The problem $\text{lex}\min\{F_1(S),\ldots,F_m(S):S\in\mathcal{S}\}$ minimizes $m$ objective functions $F_1,\ldots,F_m: \mathcal{S} \rightarrow \mathbb{R}_{0}^+$ over a set $\mathcal{S}$ of feasible solutions.
The functions are sorted in decreasing priority order, i.e.\ $F_i$ is more important than $F_{i'}$, for $i< i'$.
In a \ref{eq:lexopt} solution $S^*$, $F_1(S^*)=v_1^*=\min\{F_1(S):S\in\mathcal{S}\}$ and $F_i(S^*)=v_i^*=\min\{F_i(S):S\in \mathcal{S}, F_1(S)=v_1^*, \ldots, F_{i-1}(S)=v_{i-1}^*\}$, for $i=2,\ldots,m$.

Consider two solutions $S$ and $S'$ to the above \ref{eq:lexopt} problem. 
$S$ and $S'$ are \emph{lexicographically distinct} if there is at least one $q\in\{1,\ldots,m\}$ such that $F_q(S)\neq F_q(S')$.
Further, $S$ is \emph{lexicographically smaller} than $S'$, i.e.\ $S <_{\text{lex}} S'$ or $\vec{F}(S) <_{\text{lex}} \vec{F}(S')$, if (i) $S$ and $S'$ are lexicographically distinct and (ii) $F_q(S) < F_q(S')$, where $q$ is the smallest component in which they differ, i.e.\ $q=\min\{i:F_i(S) \neq F_i(S'), 1\leq i\leq m\}$.
$S$ is \emph{lexicographically not greater} than $S'$, i.e.\ $S \leq_{\text{lex}} S'$ or $\vec{F}(S) \leq_{\text{lex}} \vec{F}(S')$, if either $S$ and $S'$ are lexicographically equal, i.e.\ not lexicographically distinct, or $S <_{\text{lex}} S'$.
The \ref{eq:lexopt} problem $\text{lex}\min\{F_1(S),\ldots,F_m(S):S\in\mathcal{S}\}$ computes a solution $S^*$ such that $\vec{F}(S^*) \leq_{\text{lex}} \vec{F}(S)$, for all $S\in \mathcal{S}$. 

An optimal solution $S=(\vec{x},\vec{C})$ to an instance $I=(m,\mathcal{J})$ of the \ref{eq:lexopt} scheduling problem minimizes $m$ objective functions $F_1,\ldots,F_m$ lexicographically, where $F_q$ is the distinct $q$-th greatest machine completion time, for $q=1,\ldots,m$.
Lemma~\ref{Lemma:Lexicographic_Reformulation} provides an ordering of the machine completion times in a \ref{eq:lexopt} schedule and states valid inequalities.

\begin{lemma}
\label{Lemma:Lexicographic_Reformulation}
In an optimal solution to the \ref{eq:lexopt} scheduling problem: 
\begin{enumerate} 
	\item $C_i \geq C_{i+1}$, for $i=1,\ldots,m-1$, 
	\item $i\cdot C_i + \left[\sum_{q=i+1}^m C_q\right]\leq \sum_{j=1}^n p_j\leq \left[\sum_{q=1}^{i-1}C_q\right] + (m-i+1)\cdot C_i$, $\forall$ $i=1,\ldots,m$.
\end{enumerate}
\end{lemma}

\noindent
Equations (\ref{Eq:RefLex_Objective}) - (\ref{Eq:RefLex_Integrality}) formulate \ref{eq:lexopt} scheduling 
using Lemma~\ref{Lemma:Lexicographic_Reformulation}. 
\begin{subequations}
\begin{align}
\text{lex} \min_{C_i,x_{i,j}} \quad & C_1, \ldots, C_m \label{Eq:RefLex_Objective}\\
& C_i \geq C_{i+1} && M_i \in \mathcal{M} \setminus \{M_m\} \label{Eq:RefLex_MachineCompletionTimesOrder} \\
& \sum_{q=1}^{i-1} C_q + (m-i+1)\cdot C_i \geq \sum_{j=1}^n p_j && M_i \in \mathcal{M} \label{Eq:RefLex_LowerBounds} \\
& i\cdot C_i + \sum_{q=i+1}^m C_q \leq \sum_{j=1}^n p_j && M_i \in \mathcal{M} \label{Eq:RefLex_UpperBounds} \\
& C_i = \sum_{j=1}^nx_{i,j}\cdot p_j &&  M_i\in\mathcal{M} \label{Eq:RefLex_MachineCompletionTimes} \\
& \sum_{i=1}^m x_{i,j} = 1 && J_j\in\mathcal{J} \label{Eq:RefLex_JobAssignment}\\
& x_{i,j}\in\{0,1\} && J_j\in\mathcal{J}, M_i\in\mathcal{M}. \label{Eq:RefLex_Integrality}
\end{align}
\label{MILP:Lexicographic}
\end{subequations}
\vspace{-20pt}
\subsection{Perturbations}
\label{Section:Perturbation_Types}

A two-stage makespan scheduling problem is specified by an initial 
instance $I_{\init}=(m,\mathcal{J})$ and a perturbed 
instance $I_{\new}=(\hat{m},\hat{\mathcal{J}})$ of $P||C_{\max}$.
Let $\mathcal{M}$ and $\hat{\mathcal{M}}$ be the set of machines in $I_{\init}$ and $I_{\new}$, respectively.
We similarly define the sets $\mathcal{J}$ and $\hat{\mathcal{J}}$, 
denoting by $p_j$ and $\hat{p}_j$ the corresponding processing times in $I_{\init}$ and $I_{\new}$ for each job $J_j\in\mathcal{J}\cap\hat{\mathcal{J}}$.
With uncertainty realization, instance $I_{\init}$ is transformed to $I_{\new}$.
This manuscript investigates the two-stage makespan problem in the case of (i) a single perturbation, and (ii) multiple perturbations.
In the former case, the effect of uncertainty realization is one of the following perturbations: 
\begin{enumerate}
\item \emph{[Processing time reduction]} The processing time $p_j$ of job $J_j\in\mathcal{J}$ is decreased and becomes $\hat{p}_j=p_j/f_j$, for some $f_j>1$.
\item \emph{[Processing time augmentation]} The processing time $p_j$ of job $J_j\in\mathcal{J}$ is increased and becomes $\hat{p}_j=f_jp_j$, for some $f_j>1$.
\item \emph{[Job cancellation]} Job $J_j\in\mathcal{J}$ is removed, i.e.\ $\hat{\mathcal{J}}=\mathcal{J}\setminus\{J_j\}$.
\item \emph{[Job arrival]} New job $J_j\notin\mathcal{J}$ arrives, i.e.\ $\hat{\mathcal{J}}=\mathcal{J}\cup\{J_j\}$. 
\item \emph{[Machine failure]} Machine $M_i\in\mathcal{M}$ fails, i.e.\ $\hat{\mathcal{M}}=\mathcal{M}\setminus\{M_i\}$.
\item \emph{[Machine activation]} New machine $M_i\notin \mathcal{M}$ is added, i.e.\ $\hat{\mathcal{M}}=\mathcal{M}\cup\{M_i\}$. 
\end{enumerate}
These perturbations are frequently encountered in practice and investigated in the literature \cite{Kasperski2014}.
In the case of multiple perturbations, $I_{\new}$ is obtained from $I_{\init}$ by applying a series of perturbations.
Certain perturbations can be considered as equivalent.
Specifically, in some proofs: 
(i) cancelling job $J_j\in\mathcal{J}$ is identical to reducing $p_j$ to zero, i.e.\ $f_j\rightarrow\infty$,
(ii) failure of machine $M_i\in\mathcal{M}$ is equivalent to new arrivals of the jobs in $\mathcal{J}_i$, where $\mathcal{J}_i$ is the set of jobs assigned to machine $M_i$ in schedule $S_{\init}$, 
(iii) job arrivals are treated similarly to processing time augmentations.
Let $f_j$ be the \emph{perturbation factor} of job $J_j\in\mathcal{J}$. 
In our uncertainty set, $f$ is the $(k+1)$-th greatest $f_j$, $k=|\{J_{j'}\in\mathcal{J}:f_{j'}>f\}|$ is the number of \emph{unstable jobs} and $\delta=\max\{\hat{m}-m,0\}$ is the number of \emph{surplus machines} after uncertainty realization.



\section{Exact LexOpt Branch-and-Bound Algorithm (Stage 1)}
\label{Section:Lexicographic_Optimization}

This section 
introduces a \ref{eq:lexopt} branch-and-bound algorithm.
The supplementary material describes the sequential \citep{Cramer1979,Burkard1991}, weighting \citep{Sherali1982, Sherali1983}, and highest-rank objective \citep{Ogryczak1997} methods adapted to 
\ref{eq:lexopt} scheduling. 
The branch-and-bound algorithm 
uses vectorial bounds to
eliminate subtrees that cannot lexicographically dominate the incumbent, i.e.\ the best solution found thus far, by extending ideas for computing \emph{ideal points} in multiobjective optimization \citep{Erghott2006}.
In \ref{eq:lexopt} scheduling, we may derive vectorial lower and upper bounds 
by approximating a multiprocessor scheduling problem with rejections that generalizes $P||C_{\max}$. 
So, we propose packing-based algorithms for computing vectorial bounds. 
Next, we describe the branch-and-bound algorithm, our bounding approach and show their correctness. 

\begin{definition}[Vectorial Bound]
Suppose that $\vec{C}(S)=(C_1(S),\ldots,C_m(S))$ is the non-increasing vector of machine completion times in a feasible schedule $S$ of the \ref{eq:lexopt} scheduling problem.
Vector $\vec{L}=(L_1,\ldots,L_m)$ is a \emph{\bf vectorial lower bound} of $S$ if $L_i\leq C_i(S)$, for each $1\leq i\leq m$.
A \emph{\bf vectorial upper bound} $\vec{U}=(U_1,\ldots,U_m)$ of $S$ has $U_i\geq C_i(S)$, for each $1\leq i\leq m$. 
\end{definition}

\subsection{Branch-and-Bound Description}

Initially, we 
sort jobs in non-increasing processing times, i.e.\ $p_1 \geq \ldots \geq p_n$.
The search space is a tree with $n+1$ levels. 
The root node appears at level $0$.
The leaves are the set $\mathcal{S}$ of all possible $m^n$ possible schedules, i.e.\ job-to-machine assignments.
Each non-leaf node $v$ at level $\ell\in\{0,1,\ldots,n-1\}$ of the tree represents a fixed assignment of jobs $J_1,\ldots,J_{\ell}$ to the $m$ machines and jobs $J_{\ell+1},\ldots,J_n$ remain to be assigned. 
In addition, node $v$ has $m$ children corresponding to every possible assignment of job $J_{\ell+1}$ to the $m$ machines.

Denote by $\mathcal{S}(v)$ the set of all schedules in the subtree rooted at node $v$.
The branch-and-bound algorithm computes a vectorial lower bound $\vec{L}$ on the lexicographically smallest schedule $S^*\in\mathcal{S}(v)$ below node $v$.
Moreover, the primal heuristic applied in each node is longest processing time first (LPT)  
\citep{Graham1969}.
In each schedule $S$ obtained by LPT,
the branch-and-bound algorithm reorders the machines so that $C_1(S)\geq\ldots\geq C_m(S)$.
Note that this lexicographic ordering may not hold for the partial schedule of jobs $J_1,\ldots,J_{\ell}$ associated with node $v$.

Using the above components, the branch-and-bound algorithm traverses the search tree via \emph{depth-first search}.
Stack $Q$ stores the set of visited nodes that remain to be explored.
Variable $I$ stores the \emph{incumbent}, i.e.\ the current lexicographically smallest solution. 
In every step, the algorithm picks the node $u$ on top of $Q$ and explores its $m$ children. 
At each $v\in children(u)$, if LPT finds a solution $S$ such that 
$\vec{C}(S) <_{\mathrm{lex}} \vec{C}(I)$, then $I$ is updated.
If $v$ is not a leaf, Algorithm \ref{Algorithm:Vectorial_Lower_Bound} computes a vectorial lower bound $\vec{L}$ of the lexicographically best solution in $\mathcal{S}(v)$.
When $\vec{C}(I)\leq_{\mathrm{lex}}\vec{L}$, the set $\mathcal{S}(v)$ does not contain any solution lexicographically better than $I$ and the subtree rooted at $v$ is \emph{fathomed}.
Otherwise, $v$ is pushed onto stack $Q$.
Upon termination, the incumbent is optimal because every other solution has been rejected as not lexicographically smaller than the incumbent.

\subsection{Vectorial Bound Computation}
\label{Section:Vectorial_Bounds}

Next, we describe the computation of a vectorial lower bound $\vec{L} = (L_1,\ldots,L_m)$ and a vectorial upper bound $\vec{U} = (U_1,\ldots,U_m)$ at a node $v$ in the $\ell$-th level of the search tree.
Correctness proofs are provided in the supplementary material.
The algorithm performs $m$ iterations.
In iteration $i \in \{ 1, \ldots, m\}$, it calculates a lower bound $L_i$ (Algorithm \ref{Algorithm:Vectorial_Lower_Bound}) and an upper bound $U_i$ (Algorithm \ref{Algorithm:Vectorial_Upper_Bound}) on the $i$-th machine completion time using bounds $U_1,\ldots,U_{i-1}$ and $L_1,\ldots,L_{i-1}$, respectively.
Recall that $p_1\geq\ldots\geq p_n$.
W.l.o.g., each machine executes all jobs with index $\leq \ell$ before any job with index $>\ell$.
So, for each schedule in $\mathcal{S}(v)$, a unique vector $\vec{t}=(t_1,\ldots,t_m)$ specifies the machine completion times by considering only jobs $J_1,\ldots,J_{\ell}$ and ignoring the remaining ones. 
Further, no job $J_{j}$ with $\ell+1\leq j\leq n$ is executed before time $t_q$ on machine $M_q$, for $1\leq q\leq m$.

\paragraph{Vectorial lower bound component $L_i$} 
This 
computation is equivalent to constructing a \emph{pseudo-schedule} $\widetilde{S}$ where some jobs are scheduled fractionally, i.e.\ fragmented across machines.
Initially, Algorithm~\ref{Algorithm:Vectorial_Lower_Bound} fractionally assigns jobs $J_{\ell+1},\ldots,J_h$ to machines ${M_1,\ldots,M_{i-1}}$, where $h$ is the smallest index such that $\sum_{j=\ell+1}^h p_j\geq \sum_{q=1}^{i-1}(U_q-t_q)$.
For each $q = 1 \ldots i-1$, machine $M_q$ is assigned sufficiently large job pieces so that its completion time is greater than or equal to $U_q$.
Next, Algorithm \ref{Algorithm:Vectorial_Lower_Bound} fractionally assigns the remaining load $\lambda = \sum_{j=h+1}^np_j$ of jobs $J_{h+1},\ldots,J_n$ to machines $M_i,\ldots,M_m$.
This assignment minimizes the $i$-th greatest completion time in $\widetilde{S}$.
Assuming that $p_{n+1}=0$, the value $L_i$ is the maximum among $\min_{i\leq q\leq m}\{t_q\}+p_{h+1}$ and $\max_{i\leq q\leq m}\{t_q\} + \max\left\{\frac{1}{m-i+1}\left( \lambda -\sum_{q=i+1}^m(\tau-t_q)\right),0\right\}$, where $\tau=\max_{i\leq q\leq m}\{t_{q}\}$.


\begin{algorithm}[t] \nonumber
\caption{Computation of the $i$-th vectorial lower bound component}
\begin{algorithmic}[1]
\State Select the job index $\min\{h:\sum_{j=\ell+1}^{h}p_{j}\geq\sum_{q=1}^{i-1}(U_q-t_q)\}$.
\State Compute the remaining load $\lambda = \sum_{j=h+1}^np_j$.
\State Set $\tau=\max_{i\leq q\leq m}\{t_{q}\}$.
\State Return the maximum among:
\begin{itemize}
  \item $\min_{i\leq q\leq m}\{t_q\}+p_{h+1}$, and
  \item $\max_{i\leq q\leq m}\{t_q\}+\max \left\{ \frac{1}{m-i+1} \left(\lambda - \sum_{q=i+1}^m(\tau-t_q) \right),0 \right\}$.
\end{itemize}
\end{algorithmic}
\label{Algorithm:Vectorial_Lower_Bound}
\end{algorithm}


\begin{lemma}
\label{Lemma:Vectorial_Lower_Bound}
Consider a node $v$ in the $\ell$-th level of the search tree and a machine index $i\in\{1,\ldots,m\}$. Algorithm \ref{Algorithm:Vectorial_Lower_Bound} produces a value $L_i\leq C_i(S)$ for each feasible schedule $S\in\mathcal{S}(v)$ below $v$ such that $C_q(S)\leq U_q$, $\forall$ $q=1,\ldots,i-1$.
\end{lemma}

\paragraph{Vectorial upper bound component $U_i$} 
Like $L_i$, the computation of $U_i$ can be interpreted as constructing a fractional pseudo-schedule $\widetilde{S}$.
Additionally, Algorithm \ref{Algorithm:Vectorial_Upper_Bound} uses the incumbent $I$.
Schedule $\widetilde{S}$ combines the partial schedule for jobs $J_1,\ldots,J_{\ell}$ at node $v$ with a pseudo-schedule for the remaining jobs $J_{\ell+1},\ldots,J_n$ computed by Algorithm \ref{Algorithm:Vectorial_Upper_Bound}.
Initially, Algorithm \ref{Algorithm:Vectorial_Upper_Bound} assigns a total load $\sum_{q=1}^{i-1}(L_q-t_q)$ of the smallest jobs to machines $M_1,\ldots,M_{i-1}$ so that the completion time of $M_q$ becomes exactly equal to $L_q$, for $q=1,\ldots,i-1$.
That is, a piece $\tilde{p}_h$ of job $J_h$ and jobs $J_{h+1},J_{h+2},\ldots,J_n$ are assigned fractionally to machines ${M_1,\ldots,M_{i-1}}$ so that $ \tilde{p}_h + \sum_{j=h+1}^n p_j = \sum_{q=1}^{i-1}(L_q-t_q)$.
Next, Algorithm \ref{Algorithm:Vectorial_Upper_Bound} assigns the remaining load $\lambda = \sum_{j=\ell+1}^{h-1} p_j + (p_h-\tilde{p}_h)$ of jobs $J_{\ell+1},\ldots, J_h$ fractionally and uniformly to the least loaded machines among $M_i,\ldots,M_m$ as follows. 
Firstly, the partial completion times are sorted so that $t_i\leq \ldots\leq t_m$.
This sorting occurs only for computing the vectorial upper bound and does not modify any partial schedule of a node in the search tree. 
Let $\mu$ be the minimum machine index such that (i) the remaining load $\lambda$ can be fractionally scheduled to machines $M_i, \ldots, M_{\mu}$ so that they end up with a common completion time $\tau = \frac{1}{\mu-i+1} \left( \sum_{q=i}^{\mu} t_q +\lambda \right)$, and (ii) the partial completion time $t_q$ of any other machine among $M_{\mu+1},$ $\ldots,M_m$ is at least $\tau$, i.e.\ $t_{\mu+1}\geq\tau$.
Then, $U_i$ is set as the minimum among $\max\{\tau+p_{\ell},t_m\}$ and $C_i(I)$.


\begin{algorithm}[t] \nonumber
\caption{Computation of the $i$-th vectorial upper bound component}
\begin{algorithmic}[1]
\State Compute the remaining load $\lambda = \sum_{j=\ell}^np_j-\sum_{q=1}^{i-1}(L_q-t_q)$.
\State Sort machines $M_i,\ldots,M_m$ so that $t_i\leq\ldots\leq t_m$.
\State Select the machine index $\min \left \{\mu: \frac{1}{\mu-i+1} \left( \sum_{q=i}^{\mu} t_q +\lambda \right) \leq t_{\mu+1}, \;\; i\leq \mu\leq m \right\}$.
\State Return the minimum among $\max\left\{\frac{1}{\mu-i+1} \left( \sum_{q=i}^{\mu} t_q +\lambda \right) +p_{\ell}, \;\; t_m\right\}$ and $C_i(I)$.
\end{algorithmic}
\label{Algorithm:Vectorial_Upper_Bound}
\end{algorithm}

\begin{lemma}
\label{Lemma:Vectorial_Upper_Bound}
Consider a node $v$ in the $\ell$-th level of the search tree and a machine index $i\in\{1,\ldots,m\}$. Algorithm \ref{Algorithm:Vectorial_Upper_Bound} produces a value $U_i\geq C_i(S)$ for each feasible schedule $S\in\mathcal{S}(v)$ below $v$ such that $C_q(S)\geq L_q$, $\forall$ $q=1,\ldots,i-1$.
\end{lemma}

Theorem \ref{thm:bnb_correct} \mitenchange{states} the correctness of our branch-and-bound algorithm. 

\begin{theorem}
\label{thm:bnb_correct}
The branch-and-bound algorithm computes a \ref{eq:lexopt} solution.
\end{theorem}
\section{Approximate Recovery Algorithm with Binding Decisions (Stage 2)}
\label{Section:Recovery}

This section presents our recovery (reoptimization) strategy (Section \ref{Section:Recovery_Strategy}) and analyzes the price of robustness for our two-stage approach in the case of single (Section \ref{Section:Single_Perturbation}) and multiple (Section \ref{Section:Multiple_Perturbations}) perturbations.

\subsection{Recovery Algorithm Description}
\label{Section:Recovery_Strategy}

A reoptimization strategy transforms the initial schedule $S_{\init}$ to a new schedule $S_{\new}$ for the perturbed instance $I_{\new}$.
Theorem \ref{Theorem:NP-hardness} shows that optimally solving this problem is already $\mathcal{NP}$-hard.

\begin{theorem}
\label{Theorem:NP-hardness}
The makespan recovery problem 
is strongly $\mathcal{NP}$-hard, even in the case of a single perturbation.
\end{theorem}

To describe our recovery strategy, Definition~\ref{Definition:Binding_Free_Decisions} 
distinguishes between \emph{binding} and \emph{free} optimization decisions. 
Binding decisions are job assignments in $S_{\init}$ which remain valid for the perturbed instance $I_{\new}$.
Free decisions are assignments of new jobs or jobs originally
assigned to machines which failed due to uncertainty. 

\begin{definition}
\label{Definition:Binding_Free_Decisions}
Consider a makespan recovery problem instance $(I_{\init},S_{\init},I_{\new})$ with $I_{\init}=(\mathcal{M},\mathcal{J})$ and $I_{\new}=(\hat{\mathcal{M}},\hat{\mathcal{J})}$. 
\begin{itemize}
\item \textbf{Binding decisions} $\{x_{i,j}:(x_{i,j}(S_{\init})=1) \wedge(i\in \hat{\mathcal{M}} \cap 
\mathcal{M}) \wedge (j\in \hat{\mathcal{J}} \cap \mathcal{J}) \}$ are variable evaluations attainable from $S_{\init}$ in the recovery process. 
\item \textbf{Free decisions} $\{x_{i,j}: (j\in\hat{\mathcal{J}}) \wedge (\nexists i'\in \mathcal{M} \cap
\hat{\mathcal{M}} : x_{i',j}(S_{\init})=1)\}$ are variable evaluations that cannot be determined from $S_{\init}$ but are needed to recover feasibility.
\end{itemize}
\end{definition}

Our recovery strategy (Algorithm~\ref{Algorithm:Recovery}) maintains all binding decisions and makes free decisions using LPT \citep{Graham1969}. 
Theoretically,
enforcing the binding decisions exploits all relevant information in $S_{\init}$ for solving the perturbed instance $I_{\new}$, thus quantifies the benefit of staying close to $S_{\init}$. 
Practically, modifying $S_{\init}$ may incur transformation costs 
and our reoptimization algorithm mitigates this overhead.
The supplementary material presents a more flexible recovery strategy with a bounded number of binding decision modifications.



\begin{algorithm}[t] \nonumber
\caption{Recovery Strategy}
\begin{algorithmic}[1]
\State Perform all binding decisions (job assignments) with respect to schedule $S_{\init}$.
\State Schedule free (unassigned) jobs using Longest Processing Time first (LPT).
\end{algorithmic}
\label{Algorithm:Recovery}
\end{algorithm}




\subsection{Single Perturbation}
\label{Section:Single_Perturbation}

\begin{figure}[t]
    \begin{subfigure}[t]{0.5\textwidth}
    	\begin{center}
        \includegraphics{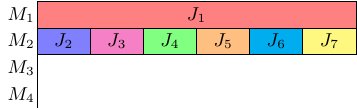}
        \end{center}
        \caption{Weakly recoverable optimal $S_{\init}$.}
        \label{Figure:Arbitrary_Initial_Schedule}
    \end{subfigure}
    \begin{subfigure}[t]{0.5\textwidth}
    	\begin{center}
        \includegraphics{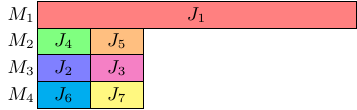}
        \end{center}
        \caption{Efficiently recoverable \ref{eq:lexopt} $S_{\init}$.}
        \label{Figure:Lexicographic_Initial_Schedule}
    \end{subfigure}
\caption{Illustration of the benefit obtained by \ref{eq:lexopt} schedules. }
\label{Figure:Lexicographic_Optimization}
\end{figure}

This section 
analyzes our two-stage approach in the case of a single perturbation.
With an arbitrary optimal initial solution, Theorem \ref{Theorem:Single_Recovery_Negative} shows that the recovery strategy results in a  non-constant price of robustness.
When the initial solution is \ref{eq:lexopt}, Theorem \ref{Theorem:Single_Recovery_Positive} provides a significantly better performance guarantee. 
Figure~\ref{Figure:Lexicographic_Optimization} illustrates a degenerate instance for deriving Theorem \ref{Theorem:Single_Recovery_Negative}, highlighting the significance of \ref{eq:lexopt} for $P||C_{\max}$ under uncertainty.

\begin{theorem}
\label{Theorem:Single_Recovery_Negative}
For the makespan recovery problem with a single perturbation, Algorithm \ref{Algorithm:Recovery} achieves an $\Omega(m)$ price of robustness with an arbitrary optimal initial schedule $S_{\init}$.
\end{theorem}

\begin{proof}
\mitenchange{Consider an instance $I_{\init}$ with $m$ machines and $n+1$ jobs, where $n=k\cdot m$, for some integer $k\in\mathbb{Z}^+$, $p_1=np$ and $p_j=p$, for $j=2,\ldots,n+1$. 
The schedule $S_{\init}$ that assigns job $J_1$ to machine $M_1$, jobs $J_2,\ldots,J_n,J_{n+1}$ to machine $M_2$ and keeps the remaining machines $M_3,\ldots,M_m$ idle, is optimal for $I_{\init}$.
Suppose that $I_{\init}$ is perturbed because job $J_1$ is cancelled and let $I_{\new}$ be the new instance (we may alternatively consider a large reduction of processing time $p_1$). 
Then, Algorithm \ref{Algorithm:Recovery} produces a schedule $S_{\rec}$ with makespan $C_{\max}(S_{\rec})=\sum_{j=1}^{n}p_j$.
However, an optimal schedule $S_{\new}$ for $I_{\new}$ has makespan $C_{\max}(S_{\new})=\frac{1}{m}\sum_{j=1}^np_j$.
Figure \ref{Figure:Lexicographic_Optimization} illustrates such two schedules, where $\frac{C_{\max}(S_{\rec})}{C_{\max}(S_{\new})}=\Omega(m)$.}
%
%
\end{proof}




Before proving 
Theorem \ref{Theorem:Single_Recovery_Positive}, we state 
Lemma~\ref{Lemma:Single_Lexicographic_Optimal_Substructure}
which relates the optimal makespan of two instances with a 
different number of machines and set of jobs. 
Denote by $C_{\max}^*(m,\mathcal{J})$ the 
optimal makespan 
for instance $(m,\mathcal{J})$.
Lemma \ref{Lemma:Single_Lexicographic_Optimal_Substructure} highlights 
the importance of a \ref{eq:lexopt} schedule: 
consider a \ref{eq:lexopt} schedule $S^*$ for $(m,\mathcal{J})$ and an arbitrary subset $\mathcal{M}'\subseteq\mathcal{M}$ of $m-1$ machines, i.e.\ $\mathcal{M}'=\mathcal{M}\setminus\{M_{\ell}\}$. 
Also, let $\mathcal{J'}$ be the subset of jobs assigned to the machines in $\mathcal{M}'$ by $S^*$.
Then, the subschedule of $S^*$ on $\mathcal{M}'$ is optimal for $(m-1,\mathcal{J}')$.

\begin{lemma}
\label{Lemma:Single_Lexicographic_Optimal_Substructure}
Consider a makespan problem instance $(m,\mathcal{J})$ and let $S$ be \ref{eq:lexopt} schedule. 
Given an arbitrary machine $M_{\ell}\in\mathcal{M}$, denote by $\mathcal{J}'$ the subset of all jobs assigned to the machines in $\mathcal{M}\setminus\{M_{\ell}\}$ by $S$. 
Then, it holds that: 
\begin{enumerate}
    \item $\max_{M_i\in\mathcal{M}\setminus\{M_{\ell}\}}\{C_i(S)\}= C_{\max}^*(m-1,\mathcal{J}')$, and
    \item $C_{\max}^*(m-1,\mathcal{J})\leq 2\cdot C_{\max}^*(m,\mathcal{J})$.
\end{enumerate}
\end{lemma}

\begin{theorem}
\label{Theorem:Single_Recovery_Positive}
For the makespan recovery problem with a single perturbation, Algorithm \ref{Algorithm:Recovery} achieves a tight price of robustness equal to 2, if $S_{\init}$ is \ref{eq:lexopt}.
\end{theorem}

\begin{proof}
The sequel proves the theorem in the case of a job reduction. 
For other perturbations described in Section \ref{Sec:Problem_Definitions}, the proof is presented in the supplementary material.
The supplementary material also shows that the obtained price of robustness is tight for every perturbation that we consider.

Consider a \ref{eq:lexopt} schedule $S_{\init}$ for instance $I_{\init}$ and suppose that the processing time of job $J_j$ decreases by $\delta\in(0,p_j]$, i.e.\ $p_j\leftarrow p_j-\delta$.
Cancelling job $J_j\in\mathcal{J}$ is equivalent to reducing $p_j$ to zero.
Suppose that machine $M_{\ell}$ executes $J_j$ in $S_{\init}$.
W.l.o.g., job $J_j$ completes last in $M_{\ell}$.
Algorithm \ref{Algorithm:Recovery} returns the recovered schedule $S_{\rec}$ which keeps the job assignments in $S_{\init}$, but decreases $p_j$ and $C_{\ell}(S_{\init})$ by $\delta$.
Let $S_{\new}$ be an optimal schedule for the perturbed instance $I_{\new}$.
We distinguish two cases depending on whether $M_{\ell}$ completes last in $S_{\rec}$, or not.

First, suppose $C_{\ell}(S_{\rec})<C_{\max}(S_{\rec})$
and let $\mathcal{J}'\subseteq \mathcal{J}$ be the subset of jobs executed by the machines in $\mathcal{M} \setminus \{M_{\ell}\}$.
Then,
{\allowdisplaybreaks
\begin{align*}
C_{\max}(S_{\rec}) = \ & C_{\max}^*(m-1,\mathcal{J}') & [\text{Lemma \ref{Lemma:Single_Lexicographic_Optimal_Substructure}.1}], \\
\leq \ & C_{\max}^*(m-1,\mathcal{J}\setminus\{J_j\}) & [\mathcal{J}'\subseteq \mathcal{J}\setminus\{J_j\}], \\
\leq \ & 2\cdot C_{\max}^*(m,\mathcal{J}\setminus\{J_j\}) & [\text{Lemma \ref{Lemma:Single_Lexicographic_Optimal_Substructure}.2}], \\
\leq \ & 2\cdot C_{\max}(S_{\new}) & [\text{Definition}].
\end{align*}
}
Subsequently, consider that $C_{\ell}(S_{\rec})=C_{\max}(S_{\rec})$, i.e.\ $C_{\max}(S_{\rec})=C_{\max}(S_{\init})-\delta$.
We claim that 
 $S_{\rec}$ is 
 optimal for $I_{\new}$.
Assume for contradiction that an optimal schedule $S_{\new}$ for $I_{\new}$ satisfies $C_{\max}(S_{\new})<C_{\max}(S_{\init})-\delta$.
By adding $\delta$ extra units of time on job $J_j$, we derive a feasible schedule $\widetilde{S}$ for $I_{\init}$ from $S_{\new}$, such that $C_{\max}(\widetilde{S})<C_{\max}(S_{\init})$.
This contradicts the optimality of $S_{\init}$ for $I_{\init}$.
\end{proof}



\subsection{Multiple Perturbations}
\label{Section:Multiple_Perturbations}

Two-stage robust optimization can be viewed as a two-player game where (i) we solve an initial instance, (ii) a malicious adversary generates perturbations, and (iii) we transform the initial solution into an efficient solution for a new instance. 
Adversarial strategies with multiple perturbations can render the initial solution weakly reoptimizable.
But \ref{eq:lexopt} can manage a bounded degree of uncertainty.
For this case, we show that Algorithm \ref{Algorithm:Recovery} produces 
solutions 
with a positive performance guarantee parameterized by the 
uncertainty set size.
Definition \ref{Def:Uncertainty_Set} describes our uncertainty set $\mathcal{U}(f,k,\delta)$ with
three parameters: 
(i) the factor $f$ indicating the boundary between \emph{stable} and \emph{unstable} job perturbations, (ii) the number $k$ of unstable jobs, and (iii) the number $\delta$ of surplus machines.
We assume that the number $k$ of unstable jobs is bounded by the number of machines $m$, i.e.\ $k < m$.

\begin{definition}
\label{Def:Uncertainty_Set}
For a makespan problem instance $(m,\mathcal{J})$ with processing times $p_1$, $\ldots$, $p_n$, the \textbf{uncertainty set} $\mathcal{U}(f,k,\delta)$ contains every instance $(\hat{m},\hat{\mathcal{J}})$ with processing times $\hat{p}_1,\ldots,\hat{p}_n$ satisfying the following properties:
\begin{itemize}
\item \textbf{Stability/instability boundary.} $\hat{\mathcal{J}}$ can be partitioned into the set $\hat{\mathcal{J}}^s$ of stable jobs and the set $\hat{\mathcal{J}}^u$ of unstable jobs, where $p_j / f \leq \hat{p}_j\leq p_j \cdot f$ $\forall$ $J_j\in\hat{\mathcal{J}}^s$.
\item \textbf{Bounded number of unstable jobs.} $|\hat{\mathcal{J}}^u|\leq k$, assuming that $k < m$.
\item \textbf{Bounded number of surplus machines.} $\max\{\hat{m}-m,0\}\leq \delta$.
\end{itemize}
\end{definition}

Suppose $C_{\max}^*(m,\mathcal{J})$ is the optimal makespan for the $P||C_{\max}$ 
instance $(m,\mathcal{J})$.
Lemma \ref{Lemma:Multiple_Lexicographic_Optimal_Substructure} (i) formalizes the optimal substructure imposed by \ref{eq:lexopt},
(ii) bounds pairwise machine completion time differences in \ref{eq:lexopt} schedules,
(iii) quantifies the 
sensitivity of the optimal makespan w.r.t.\ the number of machines, and
(iv) quantifies sensitivity of the optimal makespan  w.r.t.\ processing times.

\begin{lemma}
\label{Lemma:Multiple_Lexicographic_Optimal_Substructure}
Let $(m,\mathcal{J})$ be a makespan problem instance with a \ref{eq:lexopt} schedule $S$.
\begin{enumerate}
    \item If the subset $\mathcal{J}'\subseteq\mathcal{J}$ of jobs is executed by the subset $\mathcal{M}'\subseteq \mathcal{M}$ of machines in $S$, where $|\mathcal{M}'|=m'$, then the sub-schedule of $S$ on $\mathcal{M}'$ is optimal for $(m',\mathcal{J}')$, i.e.\ $\max_{M_i\in\mathcal{M}'}\{C_i(S)\}=C_{\max}^*(m',\mathcal{J}')$.
    \item Assuming that $M_i,M_{\ell}\in\mathcal{M}$ are two different machines such that job $J_j\in\mathcal{J}$ is assigned to $M_i$ in $S$, then $C_{\ell}(S)\geq C_i(S)-p_j$.
    \item It holds that $C_{\max}^*(m-\ell,\mathcal{J})\leq \left(1+\left\lceil \frac{\ell}{m-\ell} \right\rceil\right)\cdot C_{\max}^*(m,\mathcal{J})$ $\forall$ $\ell\in\{1,\ldots,m-1\}$.
    \item Let $(m,\hat{\mathcal{J}})$ be a makespan problem instance s.t.\ $\mathcal{J}=\hat{\mathcal{J}}$ and $\frac{1}{f}\cdot \hat{p}_j\leq p_j\leq \hat{p}_{j}$ for each $J_j$, where $p_j$ and $\hat{p}_j$ is the processing time of $J_j$ in $\mathcal{J}$ and $\hat{\mathcal{J}}$, respectively. 
Then, $\frac{1}{f}\cdot C_{\max}^*(m,\hat{\mathcal{J}}) \leq C_{\max}^*(m,\mathcal{J}) \leq C_{\max}^*(m,\hat{\mathcal{J}})$.
\end{enumerate}
\end{lemma}

\begin{table}
\small
\begin{center}
\caption{
{
Performance guarantees of our two-stage approach for different perturbations, parameterized by the (i) perturbation factor $f$, (ii) number $k < m$ of unstable jobs, and (iii) number $\delta$ of surplus machines. The term $\rho$ is the product of the performance guarantees obtained for Types 1-3.}}
\label{Table:Performance_Guarantees}
\begin{tabular}{ | lll | }
  \hline
  Type & Perturbation type 
  & Performance guarantee \\
  \hline
  Type 1 & Job cancellations, Processing time reductions & $2f\cdot(1+\lceil\frac{k}{m-k}\rceil)$ \\
  Type 2 & Processing time augmentations & $f+k$ \\
  Type 3 & Machine activations & $(1+\lceil \delta/m\rceil)$ \\
  Type 4 & Job arrivals, Machine failures & $\max\{2,\rho\}$ \\
  \hline  
\end{tabular}
\end{center}
\end{table}

Table \ref{Table:Performance_Guarantees} lists 
performance guarantees for our two-stage approach, obtained by individually analyzing each type of perturbation. 
Despite 
distinguishing the arguments for each type of perturbation, we obtain
a global price of robustness for all perturbations simultaneously by propagating the solution degradation with respect 
to the order of Table \ref{Table:Performance_Guarantees}.
Considering perturbations in this order is only for analysis purposes and does not restrict our uncertainty model.
Theoretically, \ref{eq:lexopt} is essential only for bounding the solution degradation due to job removals and processing time reductions.
But practically, the optimal substructure imposed by \ref{eq:lexopt} is beneficial in an integrated setting with all possible perturbations.
Section \ref{Section:Numerical_Results} complements the theoretical analysis with numerical experiments highlighting the 
 significance of \ref{eq:lexopt} in the recovered solution quality.
Theorem \ref{Theorem:PG} quantifies the price of robustness for our two-stage approach.

\begin{theorem}
For the two-stage robust makespan scheduling problem with $\mathcal{U}(f,k,\delta)$ uncertainty and $k<m$, our \ref{eq:lexopt}-based approach achieves a price of robustness:
\begin{equation*}
2f\cdot\left(1+\left\lceil\frac{k}{m-k}\right\rceil\right)
\cdot(f+k)\cdot\left(1+\left\lceil\frac{\delta}{m}\right\rceil\right).
\end{equation*}
\label{Theorem:PG}
\end{theorem}


\begin{proof}
Next, we analyze the recovered solution after job cancellations and processing time reductions (Type 1).
The supplementary material completes the theorem's proof for the other perturbations (Types 2-4).
Further, the supplementary material shows that the obtained price of robustness is asymptotically tight.

Processing time reductions are only recovered using binding decisions.
A job cancellation is equivalent to reducing the processing time to zero.
Given the recovered schedule $S_{\rec}$, we partition the machines $\mathcal{M}$ into the sets $\mathcal{M}^s$ of \emph{stable machines}, which are not assigned unstable jobs, and $\mathcal{M}^u$ of \emph{unstable machines}, which are assigned unstable jobs.
That is, $C_i(S_{\rec})\geq \frac{1}{f}\cdot C_i(S_{\init})$, for $M_i\in\mathcal{M}^s$, and $m^s=|\mathcal{M}^s|\geq m-k$.
Also, $\mathcal{M}^u=\mathcal{M}\setminus\mathcal{M}^s$ and $m^u=|\mathcal{M}^u|\leq k$.
Machine $M_i\in\mathcal{M}$ is \emph{critical}, if it completes last in schedule $S_{\rec}$, i.e.\ $C_i(S_{\rec})=C_{\max}(S_{\rec})$.
We distinguish two cases based on whether $\mathcal{M}^s$ contains a critical machine, or not.

\paragraph{Case 1: $\mathcal{M}^s$ contains a critical machine} 
Let $\mathcal{J}_{\new}^s\subseteq\mathcal{J}$ be the jobs assigned to machines $\mathcal{M}^s$ by $S_{\rec}$. 
Each job in $\mathcal{J}_{\new}^s$ is perturbed by a factor of at most $f$.
Let $\mathcal{J}_{\init}^s$ denote the same jobs before uncertainty realization.
Jobs in $\mathcal{J}_{\init}^s$ are executed on $\mathcal{M}^s$ in $S_{\init}$ and appear in $\mathcal{J}_{\new}^s$ with smaller processing times. 
Then,
{\allowdisplaybreaks
\begin{alignat*}{2}
C_{\max}(S_{\rec}) & = \max_{M_i\in\mathcal{M}^s}\{C_i(S_{\rec})\} \tag*{[$\mathcal{M}^s$ contains a critical machine],} \label{Eq:Multiple_Proof_Criticality} \\
& \leq \max_{M_i\in\mathcal{M}^s}\{C_i(S_{\init})\} \tag*{\text{[Processing time reduction]},} \\ 
& = C_{\max}^*(m^s,\mathcal{J}_{\init}^s) \tag*{\text{[Lemma \ref{Lemma:Multiple_Lexicographic_Optimal_Substructure}.1]},} \\ 
& \leq f\cdot C_{\max}^*(m^s,\mathcal{J}_{\new}^s) \tag*{\text{[Lemma \ref{Lemma:Multiple_Lexicographic_Optimal_Substructure}.4]},} \\ 
& \leq f\cdot C_{\max}^*(m^s,\mathcal{J}_{\new}) \tag*{\text{[$\mathcal{J}_{\new}^s\subseteq\mathcal{J}_{\new}$]},} \\ 
& = f\cdot C_{\max}^*(m-m^u,\mathcal{J}_{\new}) \tag*{[$m^s=m-m^u$],} \\ 
& \leq f\cdot\left(1+\left\lceil\frac{m^u}{m-m^u}\right\rceil\right)\cdot C_{\max}^*(m,\mathcal{J}_{\new}) \tag*{\text{[Lemma \ref{Lemma:Multiple_Lexicographic_Optimal_Substructure}.3]},} \\ 
& \leq f\cdot\left(1+\left\lceil\frac{k}{m-k}\right\rceil\right)\cdot C_{\max}(S_{\new}) \tag*{[$m^u\leq k$]}.
\end{alignat*}}
%
\paragraph{Case 2: Only $\mathcal{M}^u$ contains critical machines}
Consider an unstable critical machine $M_i\in\mathcal{M}^u$ in $S_{\rec}$, i.e.\ $C_{\max}(S_{\rec})=C_i(S_{\rec})$.
If only one job is assigned to $M_i$, then schedule $S_{\rec}$ is optimal.
Now, assume that at least two jobs are assigned to $M_i$ in $S_{\rec}$ .
Since processing times are only reduced, $C_i(S_{\rec})\leq C_i(S_{\init})$.
Because
$k<m$, there exists a machine $M_{\ell}\in \mathcal{M}^s$.
Furthermore, since $S_{\init}$ is \ref{eq:lexopt}, Lemma \ref{Lemma:Multiple_Lexicographic_Optimal_Substructure}.2 ensures that $C_{\ell}(S_{\init})\geq C_i(S_{\init})-p_j$, for each $J_j\in\mathcal{J}$ assigned to $M_i$ by $S_{\init}$.
As $S_{\init}$ contains at least two jobs, there exists a job $J_j$ assigned to $M_i$ by $S_{\init}$ such that $p_j\leq \frac{1}{2}\cdot C_i(S_{\init})$. 
Hence, $C_i(S_{\init})\leq 2\cdot C_{\ell}(S_{\init})$.
We conclude that $C_{\max}(S_{\rec})\leq 2\cdot C_{\ell}(S_{\init})$.
Because $M_{\ell}\in \mathcal{M}^s$, similarly to Case 1, we get: 
\begin{equation*}
C_{\max}(S_{\rec})\leq 2f\cdot\left(1+\left\lceil\frac{k}{m-k}\right\rceil\right)\cdot C_{\max}(S_{\new}).
\end{equation*}
%
%
\end{proof}

\section{Numerical Results}
\label{Section:Numerical_Results}

Section \ref{Section:Lexicographic_Instances_Main} describes our system specifications and the generation of benchmark $P||C_{\max}$ instances.
Section \ref{Section:Lexicographic_Experiments_Main} evaluates the \ref{eq:lexopt} branch-and-bound algorithm.
Section \ref{Section:Uncertainty_Realization_Main} presents the generation of perturbed instances, i.e.\ the effect of uncertainty realization.
Section \ref{Section:Recovery_Experiments_Main} evaluates the price of robustness of our two-stage robust scheduling approach. 

\subsection{System Specification and Benchmark Instances}
\label{Section:Lexicographic_Instances_Main}

We ran all computations on an Intel Core i7-4790 CPU 3.60GHz, 15.6 GB RAM machine running Ubuntu 14.04 64-bit.
Using Python 2.7.6 and Pyomo 4.4.1, 
we solve 
MILP models with CPLEX 12.6.3 and Gurobi 6.5.2.
The source code and test cases are available on GitHub \citep{source_code}.
We have randomly generated 
$P||C_{\max}$ instances.
\emph{Well-formed instances} admit an optimal schedule close to a \emph{perfect solution} where all machine completion times are equal, i.e.\ $C_i=C_{i'}$ for $M_i,M_{i'}\in\mathcal{M}$. 
\emph{Degenerate instances} have a less-balanced optimal schedule.
This section investigates well-formed instances and we complete the analysis with degenerate instances in the supplementary material.


Well-formed instances depend on 3 parameters: (i) the number $m$ of machines, (ii) the number $n$ of jobs, and (iii) a processing time seed $q$.
Using the parameter values in Table \ref{Table:Well-Formed_Instances_Main}, we generated \emph{moderate}, \emph{intermediate} and \emph{hard} well-formed instances.
For each combination of $m$, $n$ and $q$, we generate 3 instances based on 3 different distributions of processing times:
\emph{uniform distribution} 
$p_j \sim \mathcal{U}(\{1,\ldots,q\})$, 
\emph{normal distribution} 
$p_j \sim \mathcal{N}(q,q/3)$ and 
a \emph{symmetric of normal distribution} 
s.t.\ $p \sim \mathcal{N}(q,q/3)$ and 
$p_j=q-p$ if $p\in[0,q]$, or $p_j=2q-(p-q)$ if $p_j\in(q,2q]$.
Each processing time is rounded to the nearest integer.
Further, (symmetric) normal processing times outside $[0,2q]$ are rounded to the nearest of $0$ and $2q$.


\begin{table}[t]
\scriptsize
\begin{center}
\caption{Well-formed Instances}
\label{Table:Well-Formed_Instances_Main}
\begin{tabular}{ |c|c|c|c| } 
\hline
\textbf{Instances} & $\mathbf{m}$ & $\mathbf{n}$ & $\mathbf{q}$ \\
\hline
Moderate & $3,4,5,6$ & $20,30,40,50$ & $100,1000$ \\
\hline
Intermediate & $10,12,14,16$ & $100,200,300,400$ & $10000,100000$ \\
\hline
Hard & $10,15,20,25$ & $200,300,400,500$ & $10000,100000$ \\
\hline
\end{tabular}
\end{center}


\end{table}

\subsection{\ref{eq:lexopt} Branch-and-Bound Algorithm Evaluation}
\label{Section:Lexicographic_Experiments_Main}

This section 
numerically evaluates our branch-and-bound algorithm.
We complement our evaluation with
the sequential, highest-rank objective, and weighting methods adapted to \ref{eq:lexopt} scheduling.
Our termination criteria for each MILP solving is: 
(i) $10^3$ CPU seconds time limit, and (ii) $10^{-4}$ relative error tolerance, where the relative gap $(Ub-Lb)/Ub$ is computed using the best-found incumbent $Ub$ and the lower bound $Lb$.


The \emph{sequential method} solves $m$ MILP instances with repeated CPLEX calls, 
using the CPLEX reoptimize feature in each call to
exploit information obtained from previous calls. 
If 
$10^3$ CPU seconds in total elapse, then the method terminates when the ongoing MILP run is completed.
The \emph{highest-rank objective method} solves the \ref{eq:lexopt} scheduling problem using the CPLEX solution pool feature. 
Initially, the standard MILP model for $P||C_{\max}$ is solved. 
Next, the tree exploration goes on and populates a pool with 2000 solutions. 
The \emph{weighting method} computes a schedule $S$ 
of minimal weighted value $W(S)=\sum_{i=1}^m B^{m-i}\cdot C_i(S)$. 
We set $B = 2$ and solve the resulting MILP models with CPLEX and Gurobi. 

\begin{figure}[!ht]
\centering
    \begin{subfigure}[h]{\textwidth}
        \centering
        \includegraphics{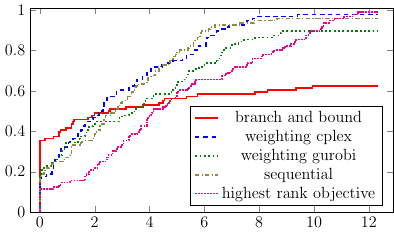}
        \includegraphics{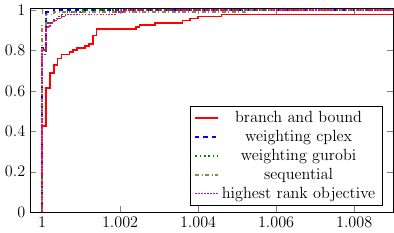}
        \caption{Moderate instances: time (s) on $\log_2$ scale (left), upper bounds on $[1,1.009]$ (right).}
    \end{subfigure}
    \begin{subfigure}[h]{\textwidth}
        \centering
        \includegraphics{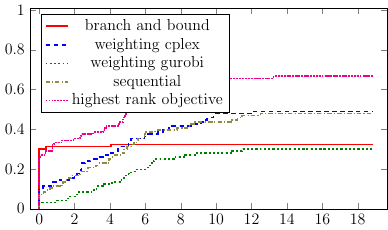}
        \includegraphics{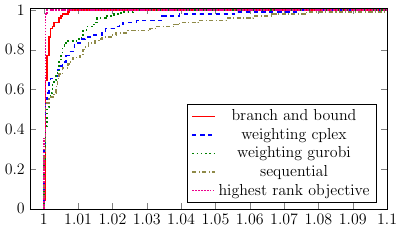}
        \caption{Intermediate instances: time (s) on $\log_2$ scale (left), upper bounds on $[1,1.1]$ (right).}
    \end{subfigure}
    \begin{subfigure}[h]{\textwidth}
        \centering
        \includegraphics{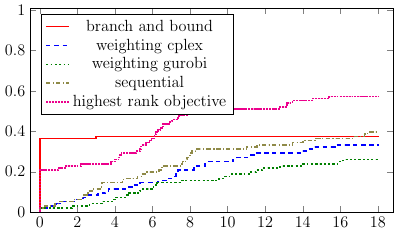}
        \includegraphics{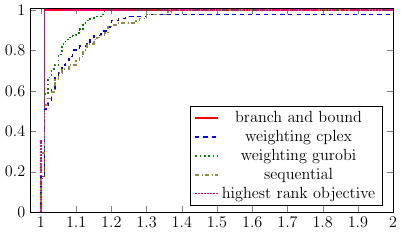}
        \caption{Hard instances: time (s) on $\log_2$ scale (left), upper bounds on $[1,2]$ (right).}
    \end{subfigure}
\caption{Performance profiles for the \emph{well-formed test set} with $10^3$ s timeout.}
\label{Figure:Wellformed_Main}
\end{figure}

Figure \ref{Figure:Wellformed_Main} plots \emph{performance profiles} of \ref{eq:lexopt} methods on well-formed instances in order to compare running times and computed solutions \citep{Dolan2002}.
The sequential and weighting methods perform similarly in terms of running time and number of solved instances.
But, the sequential method produces worse feasible solutions 
since lower-ranked objectives are not optimized in the case of a 
timeout.
The highest-rank objective method has worse running times than the sequential and weighting methods on moderate instances, because 
populating the solution pool 
is a significant part of the overall running time.
However, the highest-rank objective method attains significantly better running times than the sequential and weighting methods on intermediate and hard test cases, because populating the solution pool is a small fraction of the overall running time.
The highest-rank objective method does not prove lexicographic optimality as it only generates 2000 candidate solutions.
Nonetheless, it produces the best heuristic results for most test cases.

Figure \ref{Figure:Wellformed_Main} shows that our branch-and-bound algorithm proves global optimality faster than the other methods when it converges, i.e.\ when it does not timeout. 
Note that the branch-and-bound algorithm converges for $> 60\%$ of the moderate test cases and $> 30\%$ of the intermediate and hard instances.
For intermediate and hard instances, the branch-and-bound algorithm consistently produces good heuristic solutions, i.e.\ better solutions than the sequential and weighting methods.

\subsection{Generation of Perturbed Instances}
\label{Section:Uncertainty_Realization_Main}

Recall that an instance of the makespan recovery problem is specified by: (i) an initial makespan problem instance $I_{\init}$, (ii) an initial solution $S_{\init}$ to $I_{\init}$, and (iii) a perturbed instance $I_{\new}$. 
We generate the initial $P||C_{\max}$ instances according to Section \ref{Section:Lexicographic_Instances_Main}.
For each instance $I_{\init}$, we generate a set $\mathcal{S}(I_{\init})$ of 50 schedules using the CPLEX solution pool feature. 
In general, two different schedules $S_1,S_2\in\mathcal{S}(I_{\init})$ have weighted values $W(S_1)\neq W(S_2)$, where $W(S)=\sum_{i=1}^m B^{m-i}\cdot C_i(S)$.
For each makespan problem instance $I_{\init}$, we construct a perturbed instance $I_{\new}$ by generating random disturbances.
A \emph{job disturbance} is (i) a new job arrival, (ii) a job cancellation, (iii) a processing time augmentation, or (iv) a processing time reduction.
A \emph{machine disturbance} is (i) a new machine activation, or (ii) a machine failure.
We randomly generate $d_n = \lceil 0.2\cdot n \rceil$ job disturbances and $d_m = \lceil 0.2\cdot m \rceil$ machine perturbations.
The type of each \emph{job disturbance} is chosen uniformly at random among the four options (i) - (iv).
A new job $J_j\in\hat{\mathcal{J}}$ 
has processing time $\hat{p}_j\sim\mathcal{U}(\{1,\ldots,q\})$,
using the parameter $q$ 
that produced the original instance $I_{\init}$.
A job cancellation deletes 
an existing job chosen uniformly at random.
To increase or decrease the processing time of job $J_j\in\mathcal{J}$, 
we randomly select $\hat{p}_j\sim\mathcal{U}(\{p_j+1,\ldots,2\cdot q\})$ or $\hat{p}_j\sim\mathcal{U}(\{1,2,\ldots p_j-1\})$, respectively.
The type of a \emph{machine disturbance} is chosen uniformly at random among options (i)-(ii).
A new machine activation increases the number of available machines by one.
A machine cancellation deletes an existing machine chosen uniformly at random.

\subsection{Two-Stage Robust Scheduling Evaluation}
\label{Section:Recovery_Experiments_Main}

\begin{figure}[t]
\centering
    \begin{subfigure}{0.49\textwidth}
        \centering
        \includegraphics[width=\textwidth]{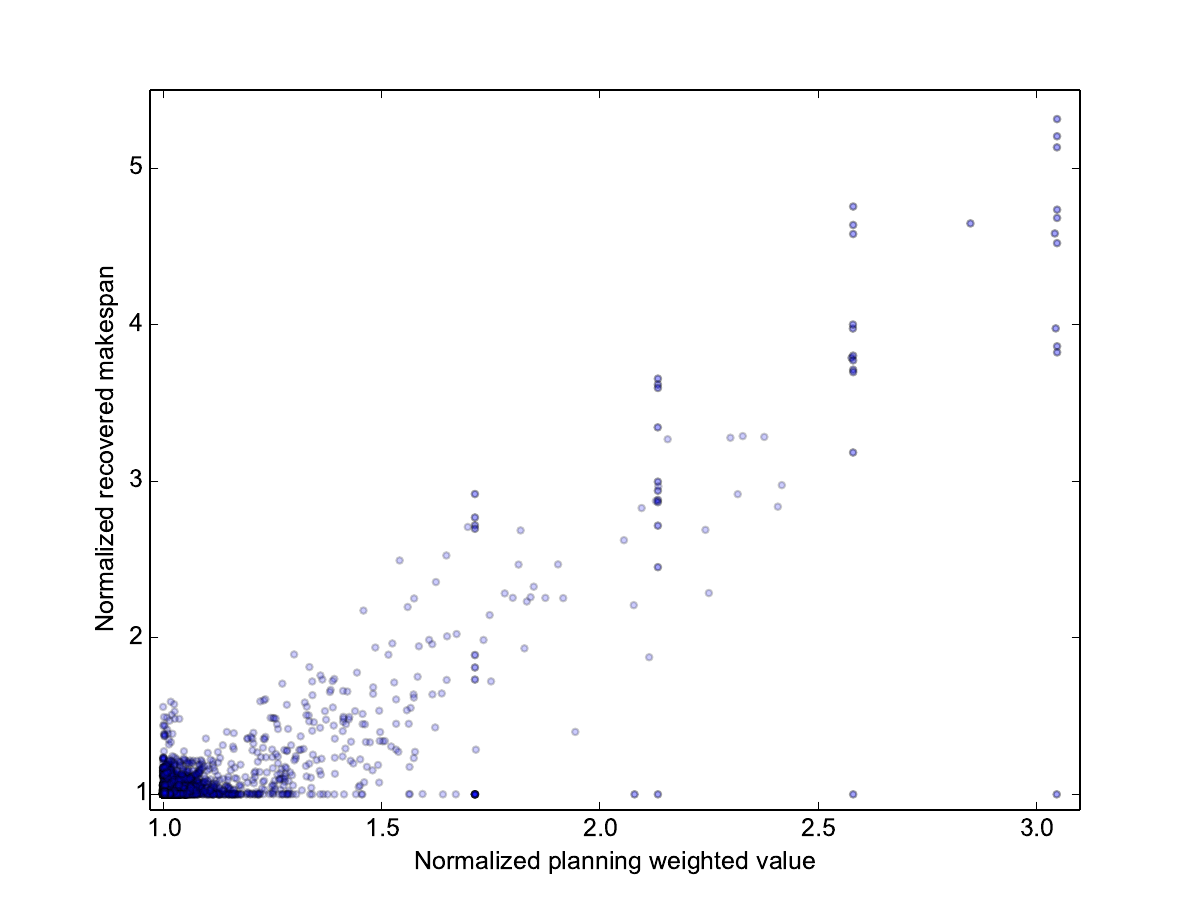}
        \caption{Binding Recovery}
        \label{Figure:Wellformed_Binding_Scatter_Plot_Main}
    \end{subfigure}
    \begin{subfigure}{0.49\textwidth}
        \centering
        \includegraphics[width=\textwidth]{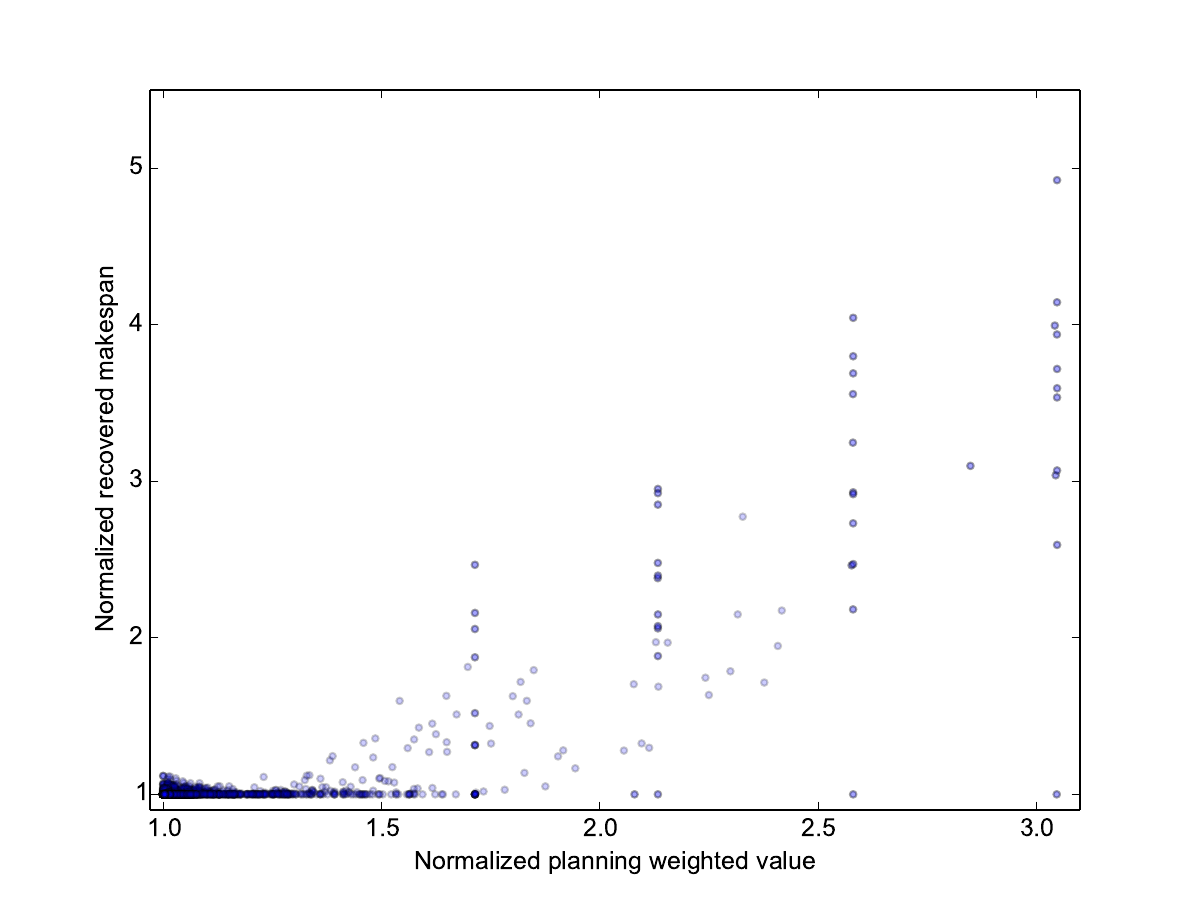}
        \caption{Flexible Recovery}
        \label{Figure:Wellformed_Flexible_Scatter_Plot_Main}
    \end{subfigure}
\caption{Well-formed instances scatter plots illustrating the recovered solution makespan with respect to the initial solution weighted value.}
\label{Figure:Wellformed_Scatter_Plots_Main}
\end{figure}

In the two-stage robust makespan scheduling problem, solution $S_{\init}$ is transformed to a feasible solution $S_{\rec}$ for instance $I_{\new}$.
Figure \ref{Figure:Wellformed_Scatter_Plots_Main} correlates the makespan of $S_{\rec}$ with the closeness of $S_{\init}$ to \ref{eq:lexopt}.
We quantify the closeness of $S_{\init}$ to \ref{eq:lexopt} using the weighted value $W(S_{\init})=\sum_{i=1}^m B^{m-i}\cdot C_i(S_{\init})$. 
The closest to \ref{eq:lexopt} the schedule $S_{\init}$ is, the lowest the value $W(S_{\init})$ we get.
For each instance $I_{\init}$, we recover every solution $S_{\init}\in\mathcal{S}(I_{\init})$ by applying our binding and flexible recovery strategies from Section \ref{Section:Recovery_Strategy}. 
For the flexible recovery strategy, we set $g=0.1n$, i.e.\ at most 10\% of the binding decisions can be modified. 
Suppose that the 
normalized weighted value of an initial solution $S_{\init}\in\mathcal{S}(I_{\init})$ is $W^N(S_{\init})=\frac{W(S_{\init})}{W^*(I_{\init})}$, where $W^*(I_{\init})$ is the best weighted value in the CPLEX solution pool for instance $I_{\init}$.
Similarly, assume that the normalized makespan of $S_{\rec}$ equal to $C^N(S_{\rec})=\frac{C_{\max}(S_{\rec})}{C_{\max}^*(I_{\new})}$, where $C_{\max}^*(I_{\new})$ is the makespan of the best recovered schedule for instance $I_{\new}$.
Figure \ref{Figure:Wellformed_Binding_Scatter_Plot_Main} 
shows that 
the makespan of solutions obtained with our binding recovery strategy tends to improve as the weighted value of the initial solution decreases.
Figure \ref{Figure:Wellformed_Flexible_Scatter_Plot_Main} 
verifies this trend for the flexible recovery strategy.
These results highlight the importance of \ref{eq:lexopt} towards efficient two-stage robust scheduling. 
Our findings also motivate scheduling under uncertainty where the planning and recovery stages are investigated 
together.

\section{Conclusion}
\label{Section:Conclusion}

Practical scheduling applications frequently require an initial, nominal schedule which is recovered after uncertainty realization.
But significantly modifying the nominal schedule might not be desirable in domains such as distributed computing \citep{yu2007adaptive} and 
timetabling \citep{Phillips2017}.  
To this end, we use exact \ref{eq:lexopt} scheduling for planning and approximate rescheduling for adaptability 
\citep{Ausiello2011,Bender2015,Chassein2016,Schieber2017}.

We provide new insights on the combinatorial structure of robust scheduling.
\ref{eq:lexopt} handles highly-symmetric mixed-integer optimization problems \citep{Balas2012, fischetti1988new,fischetti2009just,salvagnin2005dominance}, but our results also highlight \ref{eq:lexopt} benefits on scheduling under uncertainty.
By exploiting optimal substructure imposed by \ref{eq:lexopt}, we propose a two-stage robust makespan scheduling approach whose performance is substantiated with a price of robustness characterization.
Numerical results with randomly generated instances demonstrate that the closest to \ref{eq:lexopt} the initial solution is, the better the recovered solution quality we get.
Beyond scheduling, extensions to uncertain min-max partitioning problems, e.g.\ facility location and network design, with generalized cost functions are possible \citep{verschae2012power}.

Faced with the lack of strong lower bounding techniques for \ref{eq:lexopt} scheduling, we develop a new branch-and-bound algorithm, based on vectorial bounds. 
The algorithm (i) avoids iterative MILP solving of sequential methods, (ii) 
bypasses precision issues of weighting methods, and (iii) reduces the symmetry of highest-rank objective methods.
This approach is broadly relevant to \ref{eq:lexopt}.

\vspace*{0.5cm}

\noindent
\textbf{Acknowledgments} \\[2pt]
\noindent
We gratefully acknowledge support from Engineering \& Physical Sciences Research Council Research (EPSRC) [EP/M028240/1] and a Fellowship to RM [EP/P016871/1]. 


\bibliography{refs}

\newpage

\clearpage
\pagenumbering{arabic}
\renewcommand*{\thepage}{A\arabic{page}}

\title{Supplementary Material for ``Exact Lexicographic Scheduling and Approximate Rescheduling''}
\maketitle
\begin{appendices}

\paragraph{\bf Contents}
This document contains omitted parts of the manuscript \emph{Exact Lexicographic Scheduling and Approximate Rescheduling}. 
The document is a companion to the original manuscript for readers interested in complementary technicalities which have been omitted to better convey our main message, i.e.\ the importance of LexOpt in scheduling under uncertainty and relevant challenges. 
These technicalities are essential for the completeness of the presented study. 
The manuscript itself and this supplementary document 
cover the topics in a similar order.

Appendix~\ref{Appendix:LexOpt} provides omitted parts required for designing, analyzing, and evaluating the exact LexOpt methods. 
Appendix~\ref{Section:Recovery_Hardness} shows \mitenchange{that the makespan recovery problem is $\mathcal{NP}$-hard}.
Appendices~\ref{Appendix:Single_Peturbation} and \ref{Appendix:Multiple_Perturbations} 
complete the robustness analysis of our two-stage approach
in the case of 
a single and multiple perturbations, respectively. 
Appendix~\ref{Sec:Flexible_Recovery} presents a more flexible recovery strategy.
\mitenchange{Appendix~\ref{Appendix:Numerical} completes our numerical evaluation with degenerate instances.}
Finally, Appendix~\ref{Appendix:Nomenclature} provides a table with the notation used in both documents.

\section{Exact LexOpt Methods}
\label{Appendix:LexOpt}

\mitenchange{
Section~\ref{Appendix:LexOpt_Reformulation_Lemma} proves valid inequalities and
Section~\ref{Appendix:State_of_the_Art_LexOpt} adapts the sequential, weighting, and highest-rank objective methods for the \ref{eq:lexopt} scheduling problem.
Section~\ref{Appendix:Branch_and_Bound} provides a pseudo-code, 
Section~\ref{Appendix:Longest_Processing_Time_First} describes a primal heuristic,
Section~\ref{Appendix:Vectorial_Bounds} provides correctness proofs for the vectorial bounds and
Section~\ref{Appendix:Branch_and_Bound_Optimality} states a correctness proof
for our branch-and-bound algorithm.
}

\subsection{LexOpt Scheduling Reformulation Lemma}
\label{Appendix:LexOpt_Reformulation_Lemma}

\begingroup
\def\thelemma{\ref{Lemma:Lexicographic_Reformulation}}
\begin{lemma}
In an optimal solution to the \ref{eq:lexopt} scheduling problem:
\begin{enumerate} 
    \item $C_i \geq C_{i+1}$, for $i=1,\ldots,m-1$, 
    \item 
    $i\cdot C_i + \left[\sum_{q=i+1}^m C_q\right]\leq \sum_{j=1}^n p_j\leq \left[\sum_{q=1}^{i-1}C_q\right] + (m-i+1)\cdot C_i$, $\forall$ $i=1,\ldots,m$.
\end{enumerate}
\end{lemma}
\addtocounter{lemma}{-1}
\endgroup
\begin{proof}
\mitenchange{%
In any feasible schedule, the machines can be renumbered so as to satisfy the first property.}
	For the \mitenchange{second property}, observe that $\sum_{i=1}^m C_i = \sum_{j=1}^n p_j$. 
Since $C_i\geq\ldots\geq C_m$, we get that $\sum_{q=1}^{i-1} C_q + (m-i+1)\cdot C_i \geq \sum_{j=1}^n p_j$.
Similarly, given that $C_1\geq\ldots\geq C_i$, we conclude that $i\cdot C_i + \sum_{q=i+1}^m C_q\leq \sum_{j=1}^n p_j$.
\end{proof}

\subsection{State-of-the-Art \ref{eq:lexopt} Methods}
\label{Appendix:State_of_the_Art_LexOpt}




\paragraph{Sequential Method} 
\mitenchange{This method (Algorithm \ref{Algorithm:Sequential}) iteratively minimizes the objective functions $C_1,\ldots,C_n$ 
w.r.t.\ to their priority order, over the set $\mathcal{S}$ of feasible schedules \citep{Burkard1991,Cramer1979}.
Let $v_{i}^*$ be the value of $C_i$ in a \ref{eq:lexopt} solution.
The $i$-th iteration computes $v_{i}^*$ by solving MILP~(\ref{Eq:Makespan_MILP}) with the extra constraint that the first $(i-1)$ objectives should be respectively equal to $v_{1}^*,\dots,v_{i-1}^*$. 
Warm-starting iteration $i$ with the solution at iteration $(i-1)$ improves the efficiency of the method.}

\begin{algorithm}[h] \nonumber
\caption{Sequential Method}
\begin{algorithmic}[1]
\State $v_1^*=\min\{C_1:(\vec{x},\vec{C})\in\mathcal{S}\}$.
\For {$i=2,\ldots,m$}
\State $v_i^* = \min\{C_i:x\in \mathcal{S},C_1=v_1^*,\ldots, C_{i\text{-}1}= v_{i\text{-}1}^*\}$
\EndFor
\State Return the solution computed in the last iteration.
\end{algorithmic}
\label{Algorithm:Sequential}
\end{algorithm}

\paragraph{Weighting Method}
\mitenchange{
This method (Algorithm~\ref{Algorithm:Weighting}) 
minimizes a weighted sum $\sum_{i=1}^m w_i \cdot C_i$ of the objectives $C_1,\ldots,C_n$ \citep{Sherali1982}.
Typically, $w_i=B^{m-i}$ for $i=1,2,\ldots,m$, where the \emph{big-M parameter} $B>1$ is a sufficiently large constant \citep{Sherali1983}. 
Note that the highest-rank objectives are associated with the largest weights.
Further, this weighted sum can measure the distance of any solution from the \ref{eq:lexopt} solution. 
For our numerical results, we set $B=2$.}

\begin{algorithm}[h] \nonumber
\caption{Weighting Method}
\begin{algorithmic}[1]
\State Select big-M parameter $B=2$.
\For {$i=2,\ldots,m$}
\State Set machine weight $w_i=B^{m-i}$.
\EndFor
\State Solve $\min\{\sum_{i=1}^m w_i \cdot C_i: (\vec{x}, \vec{C})\in \mathcal{S}\}$. 
\end{algorithmic}
\label{Algorithm:Weighting}
\end{algorithm}


\paragraph{Highest-Rank Objective Method}
\mitenchange{This method (Algorithm \ref{Algorithm:HHR-Obj}) computes the pool $\mathcal{P}$ of all optimal solutions for the mono-objective problem $v_1^* = \min\{ C_1 : (\vec{x},\vec{C}) \in \mathcal{S}\}$ of minimizing the highest-rank objective function $C_1$, i.e.\ the makespan, and
returns the lexicographically smallest solution $\text{lex}\min\{\vec{C}(S) : \; S\in\mathcal{P}\}$ in $\mathcal{P}$ \citep{Ogryczak1997}.
A very large solution pool can be efficiently approximated with a smaller set of solutions using the CPLEX \emph{solution pool} feature.} 
In \ref{eq:lexopt}, maintaining a single solution in the pool is sufficient, if
the current solution is always replaced by a lexicographically smaller solution.
A simple greedy lexicographic comparison algorithm checks when such an update is essential.

\begin{algorithm}[h] \nonumber
\caption{Highest-Rank Objective Method}
\begin{algorithmic}[1]
\State Solve $v_1^* = \min\{C_1 : (\vec{x},\vec{C}) \in \mathcal{S}\}$.
\State Compute the solution pool $\mathcal{P} = \{ (\vec{x}, \vec{C})\in \mathcal{S} : C_1 = v_1^*\}$.
\State Return $\text{lex}\min\{\vec{C}:(\vec{x}, \vec{C})\in \mathcal{P}\}$.
\end{algorithmic}
\label{Algorithm:HHR-Obj}
\end{algorithm}

\subsection{Branch-and-Bound Algorithm Pseudocode}
\label{Appendix:Branch_and_Bound}

\begin{algorithm}[h] 
\caption{\ref{eq:lexopt} Branch-and-Bound Algorithm using Vectorial Bounds}
\label{Algorithm:Branch-and-Bound}
\begin{algorithmic}[1]
\State $Q$: empty stack
\State $r$: root node
\State $\text{push}(Q,r)$
\State $I = \{+\infty\}^m$
\While {$Q\neq\emptyset$}
\State $u=\text{top}(Q)$
\For {$v\in\text{children}(u)$}
\If {$v$ is leaf}
\State $S$: schedule of $v$
\State $I=\text{lex}\min\{I,S\}$ 
\algstore{myalg}
\end{algorithmic}
\end{algorithm}
\begin{algorithm}[h] \nonumber
\begin{algorithmic}[1]
\algrestore{myalg}
\Else
\State $S$: heuristic schedule computed via LPT
\State $I=\text{lex}\min\{I,S\}$ 
\State $\vec{L}$: vectorial lower bound of node $v$
\If {$\vec{L}\leq_{\mathrm{lex}}\vec{C}(I)$}
\State $\text{push}(Q,v)$
\EndIf
\EndIf
\EndFor
\EndWhile
\end{algorithmic}
\end{algorithm}

\subsection{Longest Processing Time First Heuristic}
\label{Appendix:Longest_Processing_Time_First}

\mitenchange{The primal heuristic applied in each node $v$ of the branch-and-bound tree is \emph{Longest Processing Time First (LPT)} (Algorithm~\ref{Algorithm:Longest_Processing_Time_First}).
%
%
LPT keeps the assignment of jobs $J_1,\ldots,J_{\ell}$, where $\ell$ is the level of node $v$, and greedily schedules jobs $J_{\ell+1},\ldots,J_n$
with the order $p_{\ell+1}\geq \ldots\geq p_n$.
In each step, LPT assigns the next job to the least-loaded machine, i.e.\ 
makes the lexicographically best decision.} 

\begin{algorithm}[h] \nonumber
\caption{Longest Processing Time First (LPT) at level $\ell$}
\begin{algorithmic}[1]
\State $\vec{t}$: Initial machine completion times
\For {$j=(\ell+1),\ldots,n$}
\State $i=\arg\min_{M_q\in\mathcal{M}}\{t_q\}$
\State $C_i \leftarrow t_i + p_j$
\EndFor
\State Sort the machines so that $C_1\geq\ldots\geq C_m$.
\end{algorithmic}
\label{Algorithm:Longest_Processing_Time_First}
\end{algorithm}

\subsection{Correctness of Vectorial Bounds}
\label{Appendix:Vectorial_Bounds}

\mitenchange{This section proves Lemmas~\ref{Lemma:Vectorial_Lower_Bound}-\ref{Lemma:Vectorial_Upper_Bound} and, thus, shows that 
Algorithms \ref{Algorithm:Vectorial_Lower_Bound}-\ref{Algorithm:Vectorial_Upper_Bound} correctly compute vectorial bounds.}



\begin{figure}[t]
    \begin{subfigure}[t]{0.65\textwidth}
        \begin{center}
        \includegraphics{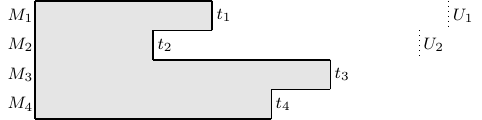}
        \end{center}
        \caption{Partial schedule associated with node $v$.}
    \end{subfigure}
    \hspace*{-1cm}
    \begin{subfigure}[t]{0.4\textwidth}
        \begin{center}
        \includegraphics{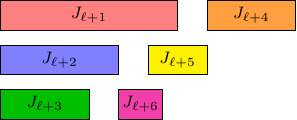}
        \end{center}
        \caption{Remaining jobs.}
    \end{subfigure}

\caption{Computing vectorial lower bound component $L_i$ at node $v$ in the $\ell$-th search tree level, by scheduling jobs $J_{\ell+1},\ldots,J_n$ in the partial schedule of $v$.
Jobs $J_{\ell+1},\ldots,J_h$ are rejected in the intervals $[t_q,U_q]$, for $q=1,\ldots,i-1$. 
$L_i$ is computed by fractionally scheduling jobs $J_{h+1},\ldots,J_m$ on machines $M_i,\ldots,M_m$ and lower bounding the completion time of machine $M_i$.} 
\label{Figure:Vectorial_Lower_Bound}
\end{figure}


\begingroup
\def\thelemma{\ref{Lemma:Vectorial_Lower_Bound}}
\begin{lemma}
Consider a node $v$ of the search tree and a machine index $i\in\{1,\ldots,m\}$. Algorithm \ref{Algorithm:Vectorial_Lower_Bound} produces a value $L_i\leq C_i(S)$ for each feasible schedule $S\in\mathcal{S}(v)$ below $v$ such that $C_q(S)\leq U_q$, $\forall$ $q=1,\ldots,i-1$.
\end{lemma}
\addtocounter{lemma}{-1}
\endgroup

\begin{proof}
	\mitenchange{%
		Schedule $S$ and pseudo-schedule $\widetilde{S}$ (of Algorithm \ref{Algorithm:Vectorial_Lower_Bound}) assign jobs $J_1,\ldots,J_{\ell}$ to the same machines and the vector $\vec{t}=(t_1,\ldots,t_m)$ specifies machine completion times w.r.t.\ these jobs.}
All remaining jobs $\mathcal{R}=\{J_{\ell+1},\ldots,J_n\}$ are scheduled differently in $\widetilde{S}$ and $S$.
In $\widetilde{S}$, 
the jobs in $\mathcal{\widetilde{R}}=\{J_{\ell+1},\ldots,J_h\}$ are fractionally assigned to machines $M_1,\ldots,M_{i-1}$ and the jobs in $\mathcal{R}\setminus\widetilde{\mathcal{R}}=\{J_{h+1},\ldots,J_n\}$ to 
$M_i,\ldots,M_m$.
Denote by $\mathcal{R}'\subseteq \mathcal{R}$ the corresponding subset of jobs assigned to machines $M_1,\ldots,M_{i-1}$, in $S$. 
That is, the jobs in $\mathcal{R}\setminus\mathcal{R}'$ are assigned to $M_i,\ldots,M_m$ in $S$.


\mitenchange{Observe that $\sum_{J_j\in\mathcal{R}'}p_j=\sum_{q=1}^{i-1}\left(C_q(S)-t_q\right)\leq \sum_{q=1}^{i-1} \left(U_q-t_q\right)$ $\leq \sum_{J_j\in\widetilde{\mathcal{R}}}p_j$, where the first equality holds by definition, the first inequality 
by the assumption 
$C_q(S)\leq U_q$, for 
$q=1,\ldots,(i-1)$, and the second inequality 
because Algorithm \ref{Algorithm:Vectorial_Lower_Bound} fits machines $M_1,\ldots,M_{i-1}$ at least up to their respective upper bounds.
Moreover, we have that $\max_{J_j\in\mathcal{R}\setminus\mathcal{R}'}\{p_j\}\geq\max_{J_j\in\mathcal{R}\setminus\widetilde{\mathcal{R}}}\{p_j\}=p_{h+1}$.
Otherwise, $\max_{J_j\in\mathcal{R}\setminus\mathcal{R}'}\{p_j\}<p_{h+1}$, which implies that
$\mathcal{R}'$ contains all jobs $J_{\ell+1},\ldots,J_{h+1}$. 
Hence, 
$\sum_{J_j\in\mathcal{R}'}p_j\geq \sum_{j=\ell+1}^{h+1}p_j> \sum_{J_j\in\widetilde{\mathcal{R}}}p_j$, 
i.e.\ a contradiction.}

\mitenchange{Since $S$ assigns a job of processing time $\max_{J_j\in\mathcal{R}\setminus\mathcal{R}'}\{p_j\}\geq p_{h+1}$ to a machine in $M_i,\ldots,M_m$,
$C_i(S)\geq \min_{i\leq q\leq m}\{t_q\}+\max_{J_j\in\mathcal{R}\setminus\mathcal{R}'}\{p_j\}\geq \min_{i\leq q\leq m}\{t_i\}+p_{h+1}$.
Clearly, $C_i(S)\geq \max_{i\leq q\leq m}\{t_i\}$.
Further, using a standard packing argument and the fact that $\sum_{J_j\in\mathcal{R}\setminus\mathcal{R}'}p_j\geq\sum_{J_j\in\mathcal{R}\setminus\widetilde{\mathcal{R}}}p_j$, if the quantity $\Lambda=\sum_{j=h+1}^n p_j-\sum_{q=i}^m(\tau-t_q)$ is positive, where $\tau=\max_{i\leq q\leq m}\{t_q\}$, then   
$C_i(S)\geq\max_{i\leq q\leq m}\{t_q\}+\frac{\Lambda}{m-i+1}$.
We conclude that $L_i \leq C_i(S)$.}
\end{proof}




\begingroup
\def\thelemma{\ref{Lemma:Vectorial_Upper_Bound}}
\begin{lemma}
Consider a node $v$ of the search tree and a machine index $i\in\{1,\ldots,m\}$. Algorithm \ref{Algorithm:Vectorial_Upper_Bound} produces a value $U_i\geq C_i(S)$ for each feasible schedule $S\in\mathcal{S}(v)$ below $v$ such that $C_q(S)\geq L_q$, $\forall$ $q=1,\ldots,i-1$.
\end{lemma}
\addtocounter{lemma}{-1}
\endgroup
\begin{proof}
\mitenchange{Recall that jobs $J_1,\ldots,J_{\ell}$ are identically assigned in schedule $S$ and pseudo-schedule $\tilde{S}$ of Algorithm \ref{Algorithm:Vectorial_Upper_Bound}. 
Moreover, a vector $\vec{t}=(t_1,\ldots,t_m)$ specifies the machine completion times of $S$ and $\tilde{S}$ w.r.t.\ these jobs.}
Let $\mathcal{R}=\{J_{\ell+1},\ldots,J_n\}$ be the set of remaining jobs.
\mitenchange{Denote by $\widetilde{\mathcal{R}}=\{J_{\ell+1},\ldots,J_h\}\subseteq\mathcal{R}$ and $\mathcal{R}'\subseteq \mathcal{R}$ the subset of jobs assigned to machines $M_i,\ldots,M_m$ in $\widetilde{S}$ and $S$, respectively.}
By arguing similarly to the proof of Lemma \ref{Lemma:Vectorial_Lower_Bound}, 
$\sum_{J_j\in\mathcal{R}\setminus\widetilde{\mathcal{R}}}p_j\geq \sum_{J_j\in\mathcal{R}\setminus\mathcal{R}'}p_j$. 
In addition, $\max_{J_j\in\mathcal{R}\setminus\widetilde{\mathcal{R}}}\{p_j\}\geq \max_{J_j\in\mathcal{R}\setminus\mathcal{R}'}\{p_j\}$. 


\mitenchange{The total load of jobs $J_{\ell+1},\ldots,J_n$ assigned to machines $M_i,\ldots,M_m$ in $S$ is clearly $\sum_{J_j\in\mathcal{R}\setminus\mathcal{R'}}p_j \leq\lambda=\sum_{j=\ell+1}^np_j-\sum_{q=1}^{i-1}(L_q-t_q)$.} 
To compute $U_i$, Algorithm \ref{Algorithm:Vectorial_Upper_Bound} assigns 
part of $\lambda$ fractionally and uniformly to the least loaded machines among $M_i,\ldots,M_m$. 
In particular, it sorts these machines so that $t_i\leq\ldots\leq t_m$ and assigns $\lambda$ units of processing time to machines $M_i,\ldots,M_{\mu}$ so that they end up having the same completion time $\tau=\frac{1}{\mu-i+1}\left(\sum_{q=1}^{\mu}t_q + \lambda\right)$.
\mitenchange{Using a simple packing argument and the fact that $\max_{J_j\in\mathcal{R}}\{p_j\}=p_{\ell}$, we get $C_i(S)\leq\max\{\tau+p_{\ell},t_m\}$.}
\end{proof}

\subsection{Optimality of Branch-and-Bound Algorithm}
\label{Appendix:Branch_and_Bound_Optimality}


\begingroup
\def\thetheorem{\ref{thm:bnb_correct}}
\begin{theorem}
The branch-and-bound method computes a \ref{eq:lexopt} solution.
\end{theorem}
\addtocounter{theorem}{-1}
\endgroup
\begin{proof}
Consider a tree node $v$. 
Let $\vec{L}=(L_1,\ldots,L_m)$ and $I$ be the computed vectorial lower bound and the incumbent, when branch-and-bound Algorithm~\ref{Algorithm:Branch-and-Bound} explores $v$.
We show the invariant that if node $v$ is pruned, then $\vec{C}(S)\geq_{\text{lex}}\vec{C}(I)$, for every schedule $S\in\mathcal{S}(v)$.
Node $v$ is pruned when $L\geq_{\text{lex}}\vec{C}(I)$, i.e.\ (i) $L_1>C_1(I)$, (ii) $L_q=C_q(I)$ $\forall$ $q=1,\ldots,i-1$ and $L_i>C_i(I)$, for some $i\in\{2,\ldots,m-1\}$, or (iii) $L_i=C_i(I)$ $\forall$ $i=1,\ldots,m$.
In case (i), because $C_1(S)\geq L_1$, it holds that $\vec{C}(S)>_{\text{lex}}\vec{C}(I)$ $\forall$ $S\in\mathcal{S}(v)$.
In case (ii), either $C_1(S)> L_1$, or $C_1(S)=L_1$ $\forall$ $S\in\mathcal{S}(v)$.
Let $\mathcal{S}_1(v)\subseteq\mathcal{S}(v)$ be the subset of schedules satisfying $C_1(S)=L_1=C_1(I)$.
Algorithm \ref{Algorithm:Vectorial_Upper_Bound} computes $U_1=C_1(I)$.
By Lemma \ref{Lemma:Vectorial_Lower_Bound}, either $C_2(S)>L_2$, or $C_2(S)=L_2$, for each $S\in\mathcal{S}_1(v)$.
Let $\mathcal{S}_2(v)\subseteq\mathcal{S}_1(v)$ be the subset of schedules with $C_2(S)=L_2$.
We define similarly all sets $\mathcal{S}_1(v),\ldots,\mathcal{S}_{i-1}(v)$.
By Lemma \ref{Lemma:Vectorial_Lower_Bound}, for any schedule in $\mathcal{S}_{i-1}(v)$, it holds that $C_q(S)=L_q=C_q(I)$ $\forall$ $q=1,\ldots,i-1$ and $C_i(S)\geq L_i>C_i(I)$.
Thus, for each $S\in\mathcal{S}(v)$, $\vec{C}(S)>_{\text{lex}}\vec{C}(I)$.
Finally, in case (iii), for each $S\in\mathcal{S}_{m-1}(v)$, $C_q(S)=C_q(I)$ $\forall$ $q=1,\ldots,m$ and $\vec{C}(S)=\vec{C}(I)$.
The theorem follows.
\end{proof}

\section{$\mathcal{NP}$-Hardness of Makespan Recovery Problem}
\label{Section:Recovery_Hardness}

\mitenchange{This section shows the $\mathcal{NP}$-hardness of the makespan recovery problem via a reduction from $P||C_{\max}$. 
That is, the makespan recovery problem is at least as hard as $P||C_{\max}$. 
Given a minimum makespan schedule $S_{\init}$ for an initial instance $I_{\init}$ of $P||C_{\max}$, the makespan recovery problem asks the existence of a feasible schedule with makespan $T_{\new}$ for a perturbed instance $I_{\new}$. 
Thus, the knowledge of $S_{\init}$ does not mitigate 
the computational complexity for solving $I_{\new}$.}
\begingroup
\def\thetheorem{\ref{Theorem:NP-hardness}}
\begin{theorem}
The makespan recovery problem 
is strongly $\mathcal{NP}$-hard, even in the case of a single perturbation.
\end{theorem}
\addtocounter{theorem}{-1}
\endgroup
\begin{proof}
\mitenchange{We prove the lemma for each type of perturbation of Section~\ref{Sec:Problem_Definitions} individually.
Given instance $I=(m,\mathcal{J})$ of $P||C_{\max}$ with target makespan $T$, 
we construct an instance $(I_{\init},S_{\init},I_{\new})$ of the makespan recovery problem 
with target makespan $T_{\new}$ by adding dummy jobs. 
Let $p_1,\ldots,p_n$ be the processing times of the $n$ jobs in $\mathcal{J}$.
Figure \ref{Figure:NP_Hardness_Appendix} shows a construction for each perturbation type.}



\begin{figure}[t]
    \begin{subfigure}[t]{0.5\textwidth}
        \begin{center}
        \includegraphics{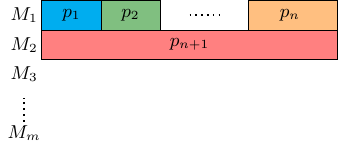}
        \end{center}
        \caption{Job cancellation, Processing time reduction.} 
        \label{Figure:NP_hardness_removal}
    \end{subfigure}
    \begin{subfigure}[t]{0.5\textwidth}
        \begin{center}
        \includegraphics{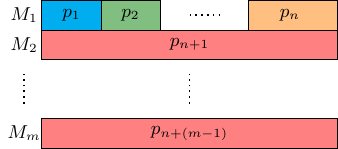}
        \end{center}
        \caption{Job arrival, Processing time augmentation, Machine failure.} 
        \label{Figure:NP_hardness_arrival}
    \end{subfigure}
    \begin{subfigure}[t]{\textwidth}
        \begin{center}
        \includegraphics{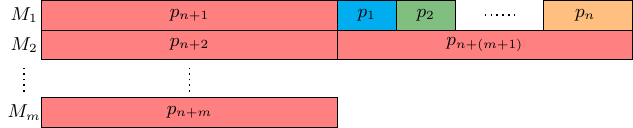}
        \end{center}
        \caption{Machine activation.} 
        \label{Figure:NP_hardness_activation}
    \end{subfigure}    
\caption{
\mitenchange{Schedule $S_{\init}$ in the $\mathcal{NP}$-hardness reduction of the makespan recovery problem from $P||C_{\max}$.  
Different perturbation types are considered individually. 
The original jobs $J_1,\ldots,J_n$ have processing times $p_1,\ldots,p_n$.
Each dummy job $J_{n+1},\ldots,J_{n+(m+1)}$ has processing time $\sum_{j=1}^n p_j$.}}
\label{Figure:NP_Hardness_Appendix}
\end{figure}

\paragraph{Job Removal, Processing Time Reduction} 
The initial instance $I_{\init}$ consists of $m$ machines, the $n$ original jobs and a dummy job of processing time $p_{n+1}=\sum_{j=1}^np_j$.
\mitenchange{Optimal schedule $S_{\init}$ for $I_{\init}$ assigns all jobs $J_1,\ldots,J_n$ to machine $M_1$, job $J_{n+1}$ to machine $M_2$ and leaves $M_3,\ldots,M_n$ empty.}
We obtain instance $I_{\new}$ from $I_{\init}$ by removing job $J_{n+1}$ \mitenchange{and setting $T_{\new}=T$}.
Since $I_{\new}$ consists only of the jobs in $I$, $I_{\new}$ admits a feasible schedule of makespan $T_{\new}$ iff there exists a schedule of makespan $T$ for $I$.
\mitenchange{The case of a processing time reduction can be treated similarly, 
by decreasing $p_{n+1}$ 
down to $0$ from $\sum_{j=1}^np_j$.}

\paragraph{Job Arrival, Processing Time Augmentation, Machine Failure} 
We construct an initial instance $I_{\init}$ with $m$ machines, the $n$ original jobs and $m-1$ 
dummy jobs $J_{n+1},\ldots,J_{n+m-1}$ 
of processing time $p_{\ell}=\sum_{j=1}^np_{j}$, for $\ell=n+1,\ldots,m-1$.
\mitenchange{The schedule $S_{\init}$ assigning 
jobs $J_1,\ldots,J_n$ to machine $M_1$ and a dummy job to every other machine is optimal for $I_{\init}$.} 
We perturb 
$I_{\init}$ by adding job $J_{n+m}$ of processing time $p_{n+m}=\sum_{j=1}^np_{j}$
and setting $T_{\new}=\sum_{j=1}^np_j+T$.
In 
instance $I_{\new}$, we ask the existence of a feasible schedule $S_{\new}$ of makespan $T_{\new}$.
Since $T<\sum_{j=1}^np_j$, if such a schedule exists, every pair of dummy jobs must executed by different machines.
\mitenchange{Thus, $I$ admits a schedule of makespan $T$ 
iff there exists a schedule with makespan $T_{\new}$ for $I_{\new}$. 
For a processing time augmentation and a machine removal, we 
use the same arguments, but different constructions. 
In the former case, we add a dummy job $J_{n+m}$ in $I_{\init}$ with $p_{n+m}=0$, which becomes $\sum_{j=1}^np_j$ in $I_{\new}$. 
In the latter case, 
we perturb 
$I_{\init}$ 
by removing 
$M_1$.}

\paragraph{Machine Activation} 
We construct a initial instance $I_{\init}$ with $m$ machines, all $n$ original jobs, and $m+1$ dummy jobs $J_{n+1},\ldots,J_{n+(m+1)}$ s.t.\ $p_{\ell}=\sum_{j=1}^np_j$, for $\ell=n+1,\ldots,n+(m+1)$.
The initial schedule $S_{\init}$ assigns a dummy job and all $n$ original jobs on machine $M_1$, two dummy jobs on machine $M_2$ and one dummy job on each machine $M_3,\ldots,M_m$.
Since any feasible schedule assigning at least two dummy jobs to one machine has makespan $\geq2\cdot\sum_{j=1}^np_j$, $S_{\init}$ must be optimal for $I_{\init}$.
We perturb $I_{\init}$ by adding a new machine and setting $T_{\new}=\sum_{j=1}^np_j+T$.
Because $T<\sum_{j=1}^np_j$, any feasible schedule for $I_{\new}$ of length $T_{\new}$ must assign one dummy job to every machine.
Thus, there exists a feasible schedule of makespan $T_{\new}$ for $I_{\new}$ iff $I$ admits a feasible schedule of makespan $T$. 
\end{proof}

\section{Robustness Analysis for a Single Perturbation}
\label{Appendix:Single_Peturbation}


This section completes the proofs of Lemma~\ref{Lemma:Single_Lexicographic_Optimal_Substructure} and Theorem~\ref{Theorem:Single_Recovery_Positive} 
for analyzing the price of robustness of our two-stage approach 
in the case of a single perturbation.

\begingroup
\def\thelemma{\ref{Lemma:Single_Lexicographic_Optimal_Substructure}}
\begin{lemma}
Consider a makespan problem instance $(m,\mathcal{J})$ and let $S$ be \ref{eq:lexopt} schedule. 
Given a machine $M_{\ell}\in\mathcal{M}$, denote by $\mathcal{J}'$ the subset of all jobs assigned to the machines in $\mathcal{M}\setminus\{M_{\ell}\}$ by $S$. 
Then, it holds that: 
\begin{enumerate}
    \item $\max_{M_i\in\mathcal{M}\setminus\{M_{\ell}\}}\{C_i(S)\}= C_{\max}^*(m-1,\mathcal{J}')$, and
    \item $C_{\max}^*(m-1,\mathcal{J})\leq 2\cdot C_{\max}^*(m,\mathcal{J})$.
\end{enumerate}
\end{lemma}
\addtocounter{lemma}{-1}
\endgroup

\begin{proof}
Suppose 
that $\max_{M_i\in\mathcal{M}\setminus\{M_{\ell}\}}\{C_i(S)\}>C_{\max}^*(m-1,\mathcal{J}')$.
\mitenchange{Let $S^*$ be 
a minimum makespan schedule for $(m-1,\mathcal{J}')$, i.e.\ $C_{\max}(S^*)=C_{\max}^*(m-1,\mathcal{J}')$.
By scheduling the jobs in $\mathcal{J}'$ as in $S^*$ and assigning the jobs in $\mathcal{J}\setminus\mathcal{J}'$ to $M_{\ell}$, we obtain a feasible schedule $\widetilde{S}$ for $(m,\mathcal{J})$ s.t.\ $\widetilde{S}<_{\text{lex}}S$, which is a contradiction.
Next, starting from an optimal schedule 
$S^*$ for $(m,\mathcal{J})$, 
we produce a new schedule $\widetilde{S}$ by moving all jobs of machine $M_m$ to 
machine $M_{m-1}$.
Clearly, $\widetilde{S}$ is a feasible for $(m-1,\mathcal{J})$ 
and the makespan has at most doubled w.r.t.\ $S^*$.}
Hence,
$C_{\max}^*(m-1,\mathcal{J})\leq C_{\max}(\widetilde{S})\leq 2\cdot C_{\max}(S^*)= 2\cdot
C_{\max}^*(m,\mathcal{J})$.
\end{proof}

\begingroup
\def\thetheorem{\ref{Theorem:Single_Recovery_Positive}}
\begin{theorem}[con't]
For the makespan recovery problem with a single perturbation, Algorithm \ref{Algorithm:Recovery} achieves a tight price of robustness equal to 2, if $S_{\init}$ is \ref{eq:lexopt}.
\end{theorem}
\addtocounter{theorem}{-1}
\endgroup


\begin{proof}
The proof of the theorem for a processing time reduction or a job removal is presented in the main manuscript.
Here, we proceed with the remaining perturbations of Section~\ref{Sec:Problem_Definitions} and show the tightness of our analysis.
Let $I_{\init}=(m,\mathcal{J})$ be the initial instance with a \ref{eq:lexopt} schedule $S_{\init}$.

\paragraph{Job Arrival, Processing Time Augmentation}
\mitenchange{Suppose that $I_{\init}$ is perturbed with the arrival of job $J_{n+1}$.
The recovered schedule $S_{\rec}$ maintains the assignments in $S_{\init}$ 
for $J_1,\ldots,J_n$ and assigns job $J_{n+1}$ to a least loaded machine 
$M_{\ell}=\arg\min_{M_i\in\mathcal{M}}\{C_i(S_{\init})\}$.
In an optimal schedule $S_{\new}$ for $I_{\new}$, it clearly holds that 
$C_{\max}(S_{\new})\geq \max\{C_{\ell}(S_{\new}),p_{n+1}\}$.
Consider the auxiliary schedule $\widetilde{S}$ obtained from $S_{\new}$ by removing job $J_{n+1}$.
Since $\widetilde{S}$ is feasible for $I_{\init}$},
$C_{\max}(S_{\new})\geq C_{\max}(\widetilde{S})\geq C_{\max}(S_{\init})$.
Hence,
$C_{\max}(S_{\rec})=\max\{C_{\ell}(S_{\init})+p_{n+1},C_{\max}(S_{\init})\}\leq C_{\max}(S_{\init})+p_{n+1}\leq 2\cdot
C_{\max}(S_{\new})$.
The case 
where $I_{\init}$ is perturbed by increasing $p_j$
can be handled using the same arguments and treating the extra piece of $J_j$ as a new job assigned to the machine executing $J_j$ in $S_{\init}$.  





\paragraph{Machine Activation, Machine Failure}
\mitenchange{Consider the case where $I_{\init}$ is perturbed because machine $M_m$ fails.}
Let $\mathcal{J}'$ be the subset of jobs assigned to $M_m$ in $S_{\init}$.
Clearly, 
$\sum_{J_j\in\mathcal{J}'}p_j\leq C_{\max}(S_{\init})$.
\mitenchange{The recovered schedule $S_{\rec}$ keeps the assignments in $S_{\init}$}
for the jobs in $\mathcal{J}\setminus\mathcal{J}'$ and assigns the jobs in $\mathcal{J}'$ to 
$M_1,\ldots,M_{m-1}$ using LPT.
Thus, $C_{\max}(S_{\rec})\leq C_{\max}(S_{\init})+\sum_{j\in\mathcal{J}'}p_j\leq 2\cdot C_{\max}(S_{\init})$.
Let $S_{\new}$ be an optimal schedule for $I_{\new}$.
\mitenchange{Since $I_{\new}$ has fewer machines 
than $I_{\init}$, 
$C_{\max}(S_{\init})\leq C_{\max}(S_{\new})$.
Hence, $C_{\max}(S_{\rec})\leq 2\cdot C_{\max}(S_{\new})$.
In the case where $I_{\init}$ is modified with the activation of a new machine $M_{m+1}$, 
$S_{\rec}$ has identical assignments with 
$S_{\init}$, while $M_{m+1}$ is left idle.
By Lemma \ref{Lemma:Single_Lexicographic_Optimal_Substructure}.2, $C_{\max}(S_{\rec})=C_{\max}^*(m,\mathcal{J})\leq 2\cdot C_{\max}^*(m+1,\mathcal{J})=2\cdot C_{\max}(S_{\new})$.}



\paragraph{Tightness}
\mitenchange{Consider an instance $I_{\init}=(m,\mathcal{J})$} with $n=m+1$ jobs of equal processing time $p$.
In a 
\ref{eq:lexopt} 
schedule $S_{\init}$, machine $M_1$ executes jobs $J_1$ and $J_2$, machine $M_i$ processes job $J_{i+1}$, for $i=2,\ldots,m$,
and $C_{\max}(S_{\init})=2p$.
\mitenchange{Assume that $I_{\init}$ is disturbed because} 
(i) job $J_n$ is removed, (ii) $p_n$ is decreased down to zero, or (iii) 
$M_{m+1}$ is activated.  
\mitenchange{In every case}, 
$C_{\max}(S_{\rec})=2p$. 
However, an optimal schedule $S_{\new}$ for $I_{\new}$, assigns exactly one job 
to each machine and 
$C_{\max}(S_{\new})=p$.
Next, consider an instance $I_{\init}=(m,\mathcal{J})$ with $n=1+(m-1)\cdot m$ jobs, 
where $p_1=m$ and 
$p_j=1$, for $j=2,\ldots,n$.
\mitenchange{In a \ref{eq:lexopt} schedule $S_{\init}$}, 
$J_1$ is assigned to machine $M_1$, exactly $m$ unit jobs are processed 
by machine $M_i$, for $i=1,\ldots,m$, 
and $C_{\max}(S_{\init})=m$.
\mitenchange{Suppose that $I_{\init}$ is perturbed because} 
(i) job $J_{n+1}$ with $p_{n+1}=m$ arrives, (ii) 
$p_n$ is augmented and becomes $m+1$, or (iii) machine $M_1$ fails.
\mitenchange{In each case}, 
$C_{\max}(S_{\rec})=2m$.
But, in an optimal schedule $S_{\new}$ for $I_{\new}$, a long job is assigned to the same machine with a unit job, 
and every other machine contains $m+1$ unit jobs, i.e.\
$C_{\max}(S_{\new})=m+1$.
\end{proof}
\section{Robustness Analysis for Multiple Perturbations}
\label{Appendix:Multiple_Perturbations}

This section completes the proofs of Lemma~\ref{Lemma:Multiple_Lexicographic_Optimal_Substructure} and Theorem~\ref{Theorem:PG} 
for analyzing the price of robustness of our two-stage approach 
in the case of multiple perturbations.





\begingroup
\def\thelemma{\ref{Lemma:Multiple_Lexicographic_Optimal_Substructure}}
\begin{lemma}
Let $(m,\mathcal{J})$ be a makespan problem instance with a \ref{eq:lexopt} schedule $S$.
\begin{enumerate}
    \item If the subset $\mathcal{J}'\subseteq\mathcal{J}$ of jobs is executed by the subset $\mathcal{M}'\subseteq \mathcal{M}$ of machines in $S$, where $|\mathcal{M}'|=m'$, then the sub-schedule of $S$ on $\mathcal{M}'$ is optimal for $(m',\mathcal{J}')$, i.e.\ $\max_{M_i\in\mathcal{M}'}\{C_i(S)\}=C_{\max}^*(m',\mathcal{J}')$.
    \item Assuming that $M_i,M_{\ell}\in\mathcal{M}$ are two different machines such that job $J_j\in\mathcal{J}$ is assigned to $M_i$ in $S$, then $C_{\ell}(S)\geq C_i(S)-p_j$.
    \item It holds that $C_{\max}^*(m-\ell,\mathcal{J})\leq \left(1+\left\lceil \frac{\ell}{m-\ell} \right\rceil\right)\cdot C_{\max}^*(m,\mathcal{J})$ $\forall$ $\ell\in\{1,\ldots,m-1\}$.
    \item Let $(m,\hat{\mathcal{J}})$ be a makespan problem instance s.t.\ $\mathcal{J}=\hat{\mathcal{J}}$ and $\frac{1}{f}\cdot \hat{p}_j\leq p_j\leq \hat{p}_{j}$ for each $J_j$, where $p_j$ and $\hat{p}_j$ is the processing time of $J_j$ in $\mathcal{J}$ and $\hat{\mathcal{J}}$, respectively. 
Then, $\frac{1}{f}\cdot C_{\max}^*(m,\hat{\mathcal{J}}) \leq C_{\max}^*(m,\mathcal{J}) \leq C_{\max}^*(m,\hat{\mathcal{J}})$.
\end{enumerate}
\end{lemma}
\addtocounter{lemma}{-1}
\endgroup
\begin{proof}
1. Assume for contradiction that $\max_{M_i\in\mathcal{M}'}\{C_i(S)\}>C_{\max}^*(m',\mathcal{J}')$.
Starting from an optimal schedule $S^*$ 
for $(m',\mathcal{J}')$, 
we construct a feasible schedule $\widetilde{S}$ for $(m,\mathcal{J})$ by assigning the jobs $\mathcal{J}'$ as in $S^*$ 
and the jobs $\mathcal{J}\setminus\mathcal{J}'$ according to $S$.
Then, $\widetilde{S}<_{\text{lex}}S$, which contradicts 
that $S$ is a \ref{eq:lexopt} schedule.

2.
\mitenchange{Assume for contradiction that $C_{\ell}(S)<C_i(S)-p_j$, i.e.\
$C_{\ell}(S)<\max\{C_\ell(S)+p_j,C_i(S)-p_j\}<C_i(S)$.
Consider the schedule $\widetilde{S}$ obtained from $S$ by moving $J_j$ from $M_i$ to $M_{\ell}$.}
Then, $C_i(\widetilde{S})=C_i(S)-p_j$, $C_{\ell}(\widetilde{S})=C_\ell(S)+p_j$, and $C_{i'}(\widetilde{S})=C_{i'}(S)$, for $M_{i'}\in\mathcal{M}\setminus\{M_i,M_{\ell}\}$.
\mitenchange{
That is, $\widetilde{S}<_{lex}S$, which contradicts that $S$ is \ref{eq:lexopt}.}

3. Starting from an optimal 
schedule $S^*$ for $(m,\mathcal{J})$, we produce a
schedule $\widetilde{S}$ by moving all jobs 
on machines $M_{m-\ell+1},\ldots,M_m$ to the remaining machines via round-robin.
For $i=1,\ldots,\ell$, the jobs of $M_{m-\ell+i}$ are moved to machine $M_{i\text{ mod }(m-\ell)}$, where $M_0=M_{m-\ell}$.
Machine $M_i\in\{M_1,\ldots,M_{m-\ell}\}$ receives jobs from at most $\lceil \ell/(m-\ell)\rceil$ machines.
Schedule $\widetilde{S}$ uses $m-\ell$ machines and its
makespan has increased 
by a factor at most $1+\lceil\frac{\ell}{m-\ell}\rceil$ w.r.t.\ $S^*$.
Hence,
$C_{\max}^*(m-\ell,\mathcal{J}) \leq C_{\max}(\widetilde{S})
\leq \left(1+\left\lceil\frac{\ell}{m-\ell}\right\rceil\right)\cdot C_{\max}(S^*) 
= \left(1+\left\lceil\frac{\ell}{m-\ell}\right\rceil\right)\cdot C_{\max}^*(m, \mathcal{J})$.

4. 
Starting from an optimal schedule $S^*$ for $(m,\mathcal{J})$, 
we construct 
a schedule $\hat{S}$ for $(m,\hat{\mathcal{J}})$ with identical assignments.
If machine $M_i$ executes a job of processing time $p_j$ in $S^*$, then $M_i$ executes a job of processing time $\hat{p}_j$ in $\hat{S}$.
Since $\frac{1}{f}\cdot\hat{p}_j\leq p_j\leq \hat{p}_j$, we have that $\frac{1}{f}\cdot C_i(\hat{S})\leq C_i(S^*)\leq C_i(\hat{S})$, for each machine $M_i$.
\end{proof}

\begingroup
\def\thetheorem{\ref{Theorem:PG}}
\begin{theorem}[con't]
For the two-stage robust makespan scheduling problem with $\mathcal{U}(f,k,\delta)$ uncertainty and $k<m$, our \ref{eq:lexopt}-based approach achieves a price of robustness:
\begin{equation*}
2f\cdot\left(1+\left\lceil\frac{k}{m-k}\right\rceil\right)
\cdot(f+k)\cdot\left(1+\left\lceil\frac{\delta}{m}\right\rceil\right).
\end{equation*}
\end{theorem}
\addtocounter{theorem}{-1}
\endgroup


\begin{proof}
\mitenchange{The proof for a processing time reduction or a job removal is presented in the main manuscript.
Here, we proceed with the remaining perturbations of Section~\ref{Sec:Problem_Definitions} and show the tightness. 
For analysis purposes, we consider the perturbations in the order of Table \ref{Table:Performance_Guarantees}. 
To propagate the solution degradation when analyzing each perturbation type, we consider that $C_{\max}(S_{\init})\leq\rho C_{\max}(S_{\init}^*)$.
That is, the initial schedule $S_{\init}$ to be recovered is $\rho$-approximate for $I_{\init}$, where $\rho\geq1$ is arbitrary.}

\paragraph{Job Cancellations, Processing Time Reductions (Type 1)}
Consider 
an instance $I_{\init}$ with $m$ machines, $(m-k)\cdot m$ jobs of length $f$, and $k$ jobs of length $m\cdot f$, \mitenchange{where $f,k=o(m)$}.
Optimal schedule $S_{\init}$ assigns $m$ jobs of length $f$ on each of the first $m-k$ machines, one job of length $m\cdot f$ on each of the remaining $k$ machines, and has makespan $C_{\max}(S_{\init})=m\cdot f$.
\mitenchange{We perturb $I_{\init}$ by decreasing the processing time of each job assigned to $M_2,\ldots,M_{m-k}$ down to 1, and cancelling the jobs assigned to the last $k$ machines.
The recovered schedule has makespan $C_{\max}(S_{\rec}) = m\cdot f$.}
But an optimal schedule $S_{\new}$ assigns to each machine a job of length $f$ and $m-k-1$ unit-length jobs, i.e.\ $C_{\max}(S_{\new})=(m-k-1)+f$.
\mitenchange{Figure \ref{Figure:Tightness_Reduction} illustrates this tightness example, where}
$\frac{C_{\max}(S_{\rec})}{C_{\max}(S_{\new})} 
= \frac{m\cdot f}{(m-k-1)+f}
= \frac{m\cdot f}{(m-k)\left(1+\frac{f-1}{m-k}\right)}
= O((1 + \frac{k}{m-k})f)$.

\begin{figure}[t]
    \begin{subfigure}[t]{0.32\textwidth}
        \begin{center}
        \includegraphics{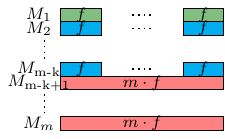}
        \end{center}
        \caption{\ref{eq:lexopt} schedule $S_{\init}$.}
        \label{Figure:Recovery_Tightness_Initial}
        \vspace*{0.5cm}
    \end{subfigure}
    \begin{subfigure}[t]{0.32\textwidth}
        \begin{center}
        \includegraphics{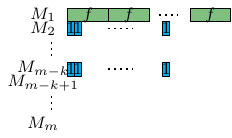}
        \end{center}
        \caption{Recovered schedule $S_{\rec}$.}
        \label{Figure:Recovery_Tightness_Recovered}
    \end{subfigure}
    \begin{subfigure}[t]{0.32\textwidth}
        \begin{center}
        \includegraphics{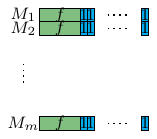}
        \end{center}
        \caption{Optimal schedule $S_{\new}$.}
        \label{Figure:Recovery_Tightness_Optimal}
        \vspace*{0.5cm}
    \end{subfigure}
\caption{
Makespan recovery instance for which the 
$O(f\cdot(1+\lceil\frac{k}{m-k}\rceil))$ factor is tight under job cancellations and processing time reductions.}
\label{Figure:Tightness_Reduction}
\end{figure}

\paragraph{Processing Time Augmentations (Type 2)}
\mitenchange{Recall that job $J_j\in\mathcal{J}$ is \emph{stable} if $\hat{p}_j\leq fp_j$ and \emph{unstable} if not.
Also, let $F=\max_{J_j\in\hat{\mathcal{J}}}\{\hat{p}_j\}$. 
Since $I_{\new}$ 
$\in\mathcal{U}(f,k,\delta)$, 
at most $k$ processing times become equal to $F$ and all remaining jobs are increased by a factor at most $f$.
Given that $S_{\rec}$ 
is identical with $S_{\init}$, 
except that some processing times are increased,
$C_i(S_{\rec})\leq f\cdot C_i(S_{\init}) + k\cdot F$ 
for each $M_i\in\mathcal{M}$.
Given an optimal schedule $S_{\init}^*$ for $I_{\init}$ and that
the processing times in $I_{\new}$ are one-to-one greater than or equal to the ones in 
$I_{\init}$,
$C_i(S_{\init})\leq\rho\cdot C_{\max}(S_{\init}^*)\leq\rho\cdot C_{\max}(S_{\new})$.
In addition, $C_{\max}(S_{\new})\geq F$.
Hence, $C_{\max}(S_{\rec})\leq (f+k)\rho\cdot C_{\max}(S_{\new})$.}

For the tightness, 
consider 
instance $I_{\init}$ (Figure \ref{Figure:Tightness_Augmentation}) with $m$ machines and $n=m^2$ unit-length jobs.
In $S_{\init}$, each machine executes $m$ jobs 
and $C_{\max}(S_{\init})=m$.
\mitenchange{After uncertainty realization, 
$k$ jobs assigned to $M_1$ get processing time $F$, 
every other job on $M_1$ gets processing time $f$,
and all other processing times remain the same in $I_{\init}$ and $I_{\new}$.}
Suppose that $F=f+m$, $F=\Theta(m)$, $f=o(m)$ and $k=o(m)$.
Schedule $S_{\rec}$ 
performs identical assignments with $S_{\init}$, i.e.\ 
$C_{\max}(S_{\rec})=k\cdot F + (m-k)\cdot f$.
\mitenchange{In an optimal schedule $S_{\new}$ for $I_{\new}$, each machine $M_i$ with $i\leq k$ processes a job of length $F$ and $k$ unit length jobs,  
while each machine $M_i$ with $i>k$ executes a job of length $f$ and $m+k$ unit length jobs, i.e.\ $C_{\max}(S_{\new})=F+m$.} 
Thus,
$\frac{C_{\max}(S_{\rec})}{C_{\max}(S_{\new})} 
= \frac{k\cdot F}{F+k} + \frac{(m-k)\cdot f}{f+m+k}
= k\cdot\frac{1}{1+\frac{k}{F}}+f\cdot\frac{1}{1+\frac{f+2k}{m-k}}=O(f+k)$.

\begin{figure}[t]
    \begin{subfigure}[t]{0.29\textwidth}
        \begin{center}
        \includegraphics{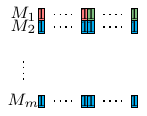}
        \end{center}
        \caption{\ref{eq:lexopt} schedule $S_{\init}$.}
        \label{Figure:Recovery_Tightness_Initial}
        \vspace*{0.5cm}
    \end{subfigure}
    \begin{subfigure}[t]{0.36\textwidth}
        \begin{center}
        \includegraphics{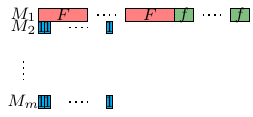}
        \end{center}
        \caption{Recovered schedule $S_{\rec}$.}
        \label{Figure:Recovery_Tightness_Recovered}
    \end{subfigure}
    \begin{subfigure}[t]{0.29\textwidth}
        \begin{center}
        \includegraphics{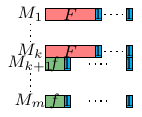}
        \end{center}
        \caption{optimal schedule $S_{\new}$.}
        \label{Figure:Recovery_Tightness_Optimal}
        \vspace*{0.5cm}
    \end{subfigure}
\caption{Makespan recovery instance for which 
the $O(f+k)$ factor is 
tight 
under processing time augmentations.}
\label{Figure:Tightness_Augmentation}
\end{figure}





\paragraph{Machine Activations (Type 3)}
\mitenchange{
Denote by $\mathcal{M}^s$ the set of available machines 
in $S_{\init}$ and by $\mathcal{M}^u$ the set of newly activated machines after uncertainty realization. 
Algorithm~\ref{Algorithm:Recovery} keeps the schedule $S_{\init}$ for the machines in $\mathcal{M}^s$ and leaves the machines in $\mathcal{M}^u$ idle, i.e.\ $C_{\max}(S_{\rec})\leq\rho\cdot C_{\max}^*(m,\mathcal{J})$.
By definition, 
$C_{\max}(S_{\new})=C_{\max}^*(m+k,\mathcal{J})$.
Hence, by Lemma \ref{Lemma:Multiple_Lexicographic_Optimal_Substructure}.3, we conclude that $C_{\max}(S_{\rec})\leq (1+\lceil k/m\rceil)\rho\cdot C_{\max}(S_{\new})$.
For the tightness of this bound, consider an instance with $m$ machines and $n=m\cdot(m+k)$ unit jobs.
In $S_{\init}$, each machine processes $m+k$ unit jobs 
and $C_{\max}(S_{\init})=m+k$.
In $S_{\rec}$, all newly activated machines are empty and $C_{\max}(S_{\rec})=m+k$.
But an optimal schedule $S_{\new}$ for $I_{\new}$ assigns exactly $m$ jobs on each machine, i.e.\ $C_{\max}(S_{\new})=m$.
That is, $\frac{C_{\max}(S_{\rec})}{C_{\max}(S_{\init})}=1+\frac{k}{m}$.
}





\paragraph{Job Arrivals, Machine Failures (Type 4)}
\mitenchange{
Consider a set 
of \emph{free jobs} arriving after uncertainty realization. 
Algorithm~\ref{Algorithm:Recovery} schedules these jobs according to LPT.
We partition the set $\mathcal{M}$ of machines
into the set $\mathcal{M}^s$ of \emph{stable machines}, not executing free jobs, and the set $\mathcal{M}^u$ of \emph{unstable machines}, assigned free jobs, in $S_{\rec}$. 
Since $\mathcal{J}_{\init}\subset\mathcal{J}_{\new}$ and $\mathcal{M}_{\init}=\mathcal{M}_{\new}$,
if there is a machine $M_i\in\mathcal{M}^s$ with $C_i(S_{\rec})=C_{\max}(S_{\rec})$, then 
$C_{\max}(S_{\rec})=C_{\max}(S_{\init})\leq \rho\cdot C_{\max}^*(I_{\init}) \leq\rho\cdot  C_{\max}(S_{\new})$.
If there is not such a machine in $\mathcal{M}^s$, we use the analysis for LPT by \citet{Graham1969}.
Since the last completing job $J_j$ 
begins at $b_j$, all machines are occupied until $b_j$ in $S_{\rec}$.
Thus, $C_{\max}(S_{\rec})=b_j+p_j\leq\frac{1}{m}\sum_{J_{j'}\in\mathcal{J}_{\new}}p_{j'}+p_j
\leq 2\cdot C_{\max}(S_{\new})$.
In both cases, $C_{\max}(S_{\rec})\leq\max\{2,\rho\}\cdot C_{\max}(S_{\new})$.
The case where a machine $M_i\in\mathcal{M}$ fails can be treated similarly by considering the jobs originally assigned to $M_i$ in $S_{\init}$ as free.
In the case where $\rho>2$, we may design a makespan problem instance such that the arrival of a new job with a tiny processing time has the effect that the recovered schedule remains $\rho$-approximate.
The tightness example for LPT implies that Algorithm \ref{Algorithm:Recovery} cannot result in a price of robustness better than 2.}
\end{proof}

\section{Flexible Recovery Strategy}
\label{Sec:Flexible_Recovery}

\mitenchange{Next, we present 
a more flexible recovery strategy (than Algorithm \ref{Algorithm:Recovery})
that modifies a bounded number of binding decisions 
\citep{CHASSEIN2018423,nasrabadi2013robust}.
To this end, we formulate the makespan recovery problem as a MILP.}
%
Let $\mathcal{J}^B=\{J_j\in\mathcal{J}_{\init}\cap\mathcal{J}_{\new}: \exists i \text{ with } x_{i,j}(S_{\text{init}}) = 1 \}$ be the binding decisions, 
i.e.\ the jobs 
appearing both in 
$I_{\init}$ and $I_{\new}$. 
Algorithm~\ref{Algorithm:Recovery} \mitenchange{keeps} the assignments in $S_{\init}$ for the \emph{binding jobs} $\mathcal{J}^B$ and greedily schedules the \emph{free jobs} $\mathcal{J}^F= \mathcal{J} \setminus \mathcal{J}^B$ \mitenchange{with} LPT \mitenchange{to produce $S_{\rec}$}.
A more flexible recovery strategy migrates a bounded number $g$ of binding jobs. %
These migrations 
produce better recovered solutions at the price of extra computational effort and higher transformation cost. 
\mitenchange{Denote by $\mathcal{J}_i^B\subseteq \mathcal{J}^B$ and $\mu_j$ the subset of binding jobs assigned to machine $M_i$ and the machine index which job $J_j\in\mathcal{J}^B$ is assigned to, respectively, in $S_{\init}$.
Our flexible recovery strategy solves MILP (\ref{Eq:Makespan_MILP}) with the additional constraint $\sum_{J_j\in\mathcal{J}^B} \sum_{M_i\in \mathcal{M}\setminus\{M_{\mu_j}\}} x_{i,j} \leq g$.}



\section{Numerical Results with Degenerate Instances}
\label{Appendix:Numerical}


This section complements our numerical results with 
degenerate instances of $P||C_{\max}$. 
Section \ref{Section:Lexicographic_Instances} describes the generation of these instances.
Sections \ref{Section:Lexicographic_Experiments}-\ref{Section:Recovery_Experiments} evaluate the \ref{eq:lexopt} branch-and-bound algorithm and the robustness of our two-stage approach using the new instances.

\subsection{Generation of Degenerate Instances}
\label{Section:Lexicographic_Instances}

\emph{Degenerate instances} have less balanced optimal solutions than \emph{well-formed instances}.
\mitenchange{
To produce degenerate instances, we sample integer processing times that can be encoded with $b$ bits.} 
Instances with small $\kappa=b/n$ values are easier to solve than instances with larger $\kappa$ values \citep{bauke2003phase}.
The \emph{phase transition} from ``easy'' to ``hard'' instances becomes sharper as $n$ increases and occurs at the threshold value $\kappa^*=\frac{\log_2m}{m-1}$. 
\mitenchange{Instances with small $\kappa$ admit exponentially many perfect solutions (where all machine completion times are equal). 
Instances with $\kappa>\kappa^*$ 
have 
less or no perfect solutions.} 
Similar phase transitions occur for other fundamental combinatorial optimization problems, e.g.\ satisfiability \citep{Mitchell1992} and the traveling salesman problem \citep{Gent1996}, where instances near the threshold value tend to be the most difficult.

\mitenchange{We derive degenerate instances by varying two parameters: (i) the number $m$ of machines and (ii) the number $n$ of jobs. 
Further, we use a processing time seed $q=2^{\lfloor\kappa(m)\cdot n\rfloor}$ with $\kappa(m)=(\log_2m)/(m-1)$.
Table \ref{Table:Degenerate_Instances} reports this information for \emph{moderate} and \emph{intermediate} degenerate instances. 
For each combination of $m$, $n$ and the corresponding $q=2^{\lfloor\kappa(m)\cdot n\rfloor}$ value, we generate 3 instances by sampling processing times from the set $\{1,\ldots,q\}$ using the uniform, normal and symmetric of normal distributions, similarly to the well-formed instances.}

\begin{table}[t]
\begin{center}
\scriptsize
\caption{Degenerate Instances}
\label{Table:Degenerate_Instances}
\begin{tabular}{ |c|c|c|c| } 
\hline
\textbf{Instances} & $\mathbf{m}$ & $\mathbf{n}$ & $\mathbf{q}$ \\
\hline
\multirow{4}*{Moderate} & 3 & $20,25,30,35$ & $2^{15},2^{19},2^{23},2^{27}$ \\ 
& 4 & $25,30,35,40$ & $2^{16},2^{20},2^{23},2^{26}$ \\ 
& 5 & $30,35,40,45$ & $2^{17},2^{20},2^{23},2^{26}$ \\ 
& 6 & $35,40,45,50$ & $2^{18},2^{20},2^{23},2^{25}$ \\ 
\hline
\multirow{4}*{Intermediate} & 10 & $40,50,60,70$ & $2^{14},2^{18},2^{22},2^{25}$\\ 
& 12 & $45,55,65,75$ & $2^{14},2^{17},2^{21},2^{24}$ \\ 
& 14 & $55,65,75,85$ & $2^{16},2^{19},2^{21},2^{24}$ \\ 
& 16 & $60,70,80,90$ & $2^{16},2^{18},2^{21},2^{24}$ \\ 
\hline
\end{tabular}
\end{center}
\end{table}

\subsection{\ref{eq:lexopt} Branch-and-Bound Algorithm Evaluation}
\label{Section:Lexicographic_Experiments}

\begin{figure}[h]
\centering
    \begin{subfigure}[h]{\textwidth}
        \centering
        \includegraphics{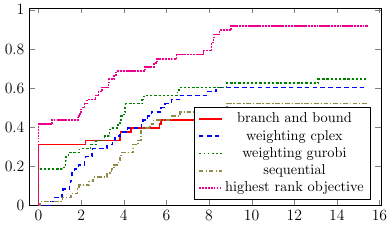}
        \includegraphics{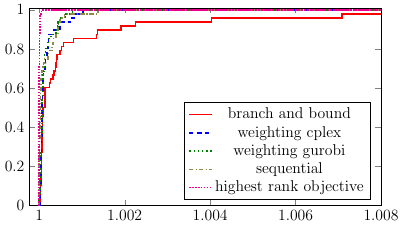}
        \caption{Moderate instances: time (s) on $\log_2$ scale (left), upper bounds on $[1,1.008]$ (right).}
        \label{Figure:Degenerate_Moderate}
    \end{subfigure}
    \begin{subfigure}[h]{\textwidth}
        \centering
        \includegraphics{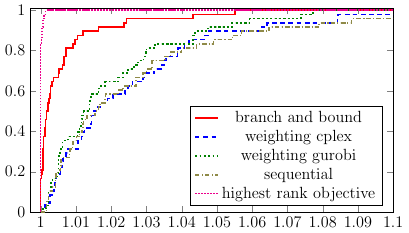}
        \caption{Intermediate instances: upper bounds on $[1,1.1]$. \emph{No solver converges for any intermediate degenerate instance within the specified time limit.}}
        \label{Figure:Degenerate_Hard}
    \end{subfigure}
\caption{Performance profiles for the \emph{degenerate test set} with $10^3$ s timeout.}
\label{Figure:Degenerate}
\end{figure}

\mitenchange{We evaluate our branch-and-bound algorithm on degenerate instances in comparison with the sequential, weighting and highest-rank objective methods. 
For MILP solving, we use (i) $10^3$ CPU seconds time limit and (ii) $10^{-4}$ error tolerance, similarly to the numerical results obtained with well-formed instances. 
Figure \ref{Figure:Degenerate} shows performance profiles for evaluating the running times and quality of computed solutions on degenerate instances.}
We observe that degenerate instances are significantly harder to solve than well-formed instances of identical size.
For instance, no solver converges for any intermediate degenerate instance, while every solver converges for $> 30\%$ of the intermediate well-formed instances.
In terms of solver comparison, we derive similar results to those obtained for  well-formed instances.
The sequential method performs similarly to the weighting method.
The highest-rank objective method produces the best heuristic results.
Our branch-and-bound method produces the second best heuristic result for intermediate degenerate instances and computes \ref{eq:lexopt} solutions quickly when it terminates. 

\begin{figure}[!ht]
\centering
    \begin{subfigure}{0.49\textwidth}
        \centering
        \includegraphics[width=\textwidth]{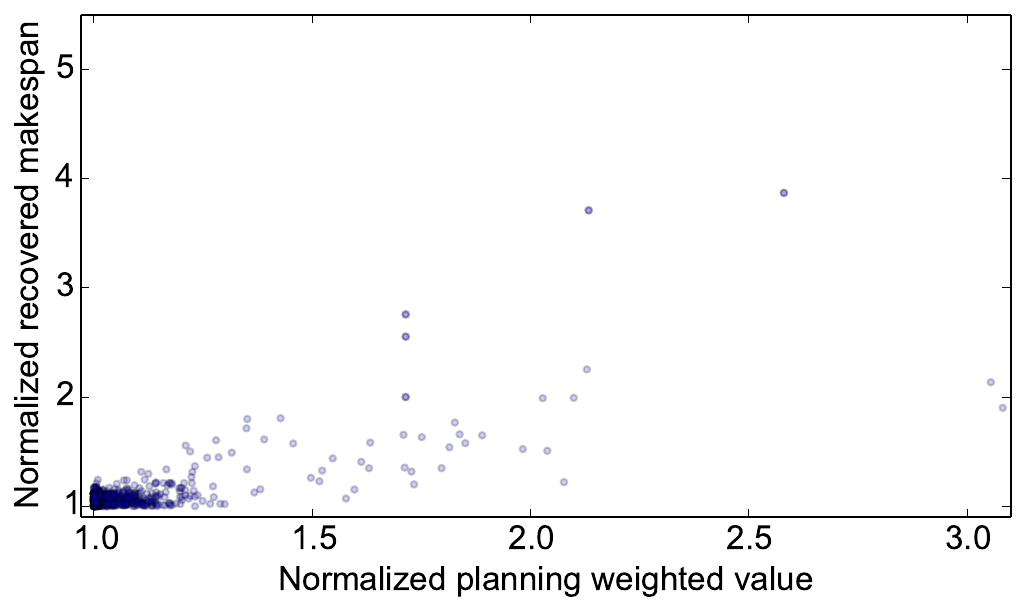}
        \caption{Binding Recovery}
        \label{Figure:Degenerate_Binding_Scatter_Plot}
    \end{subfigure}
    \begin{subfigure}{0.49\textwidth}
        \centering
        \includegraphics[width=\textwidth]{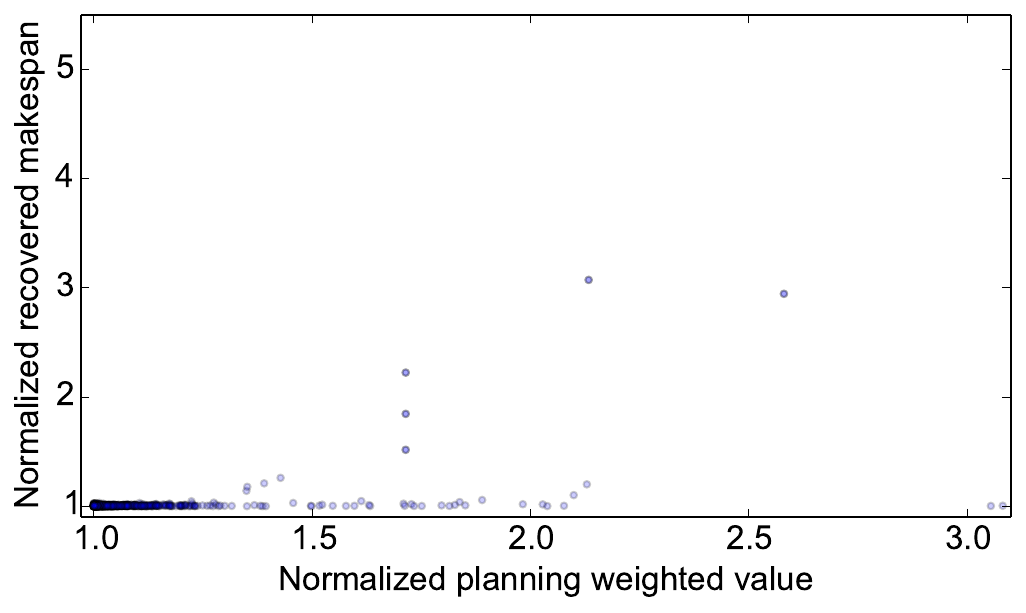}
        \caption{Flexible Recovery}
        \label{Figure:Degenerate_Flexible_Scatter_Plot}
    \end{subfigure}
\caption{Degenerate instances scatter plots illustrating the recovered solution makespan with respect to the initial solution weighted value.}
\label{Figure:Degenerate_Scatter_Plots}
\end{figure}

\subsection{Two-Stage Robustness Assessment}
\label{Section:Recovery_Experiments}

\mitenchange{Next, we investigate the impact of \ref{eq:lexopt} to the quality of the recovered solution for degenerate instances. 
For each degenerate instance $I_{\init}$ derived according to Section~\ref{Section:Lexicographic_Instances}, 
we compute 50 diverse initial solutions using the CPLEX solution pool feature. 
To quantify the closeness of an initial solution to \ref{eq:lexopt}, we use the weighted value $W(S)=\sum_{i=1}^m B^{m-i}\cdot C_i(S)$.
Further, we obtain a perturbed instance $I_{\new}$ by generating random disturbances, similarly to well-formed instances. 
Then, we fix every initial solution by applying our binding and flexible recovery strategies in Section \ref{Section:Recovery_Strategy}. 
As in the case of well-formed instances, Figures 
\ref{Figure:Degenerate_Binding_Scatter_Plot} and \ref{Figure:Degenerate_Flexible_Scatter_Plot} 
plot the normalized makespan $C^N(S_{\rec})$ obtained by our recovery strategies 
and the normalized initial solution weighted value $W^N(S_{\init})$  
for every recovered solution.
Clearly, the recovered solution improves if the initial solution weighted value decreases.
Moreover, flexibility enables more efficient recovery.
Interestingly, degenerate instances are recovered more efficiently than well-formed ones.}

\section{Table of Notation}
\label{Appendix:Nomenclature}

\begin{center}
\footnotesize
\begin{longtable}{l l}
\caption{Nomenclature}
\label{Table:Nomenclature}\\

\toprule
Name & Description \\


\midrule
\multicolumn{2}{l}{\bf Makespan scheduling problem} \\
$I=(m,\mathcal{J})$ & Instance \\
$i,q,\mu$ & Machine indices ($q,\mu$ typically used as auxiliary machine indices) \\
$j,h,\ell$ & Job indices ($h,\ell$ typically used as auxiliary job indices) \\
$m$ & Number of machines \\
$M_i \in \mathcal{M}$ & Machine $M_i$ 
in the set $\mathcal{M} = \{M_1,\ldots,M_m\}$ of all machines \\
$n$ & Number of jobs \\
$J_j \in \mathcal{J}$ & Job $J_j$ in the set $\mathcal{J} \in \{J_1,\ldots,J_n\}$ of all jobs \\ 
$p_j$ & 
Processing time of job $J_j$ \\

$C_{\max}$ & Makespan \\
$C_i$ & 
Completion time of machine $M_i$ \\ 
$x_{i,j}$ & Binary variable indicating an assignment of job $J_j$ to machine $M_i$ \\

$S=(\vec{y},\vec{C}),S',\widetilde{S}$ & Schedules ($S',\widetilde{S}$ typically used as auxiliary schedules) \\
$S^*$ & \ref{eq:lexopt} schedule \\
$\mathcal{S}$ & Set of all feasible schedules \\

\midrule
\multicolumn{2}{l}{\bf LexOpt scheduling problem} \\
$\leq_{\text{lex}}$ & Operator for lexicographic comparison \\
$F_i$ & 
$i$-th greatest completion time, i.e.\ $i$-th objective function \\
$v_i^*$ & Value of $F_i$ in a \ref{eq:lexopt} schedule $S^*$ \\ 
$\mathcal{T}_q$ & Set of tuples $(i_1,\ldots,i_q)$ with $q$ pairwise disjoint machine indices \\
$w_i$ & Weight of objective function $F_i$ (weighting method) \\
$\mathcal{P}$ & Solution pool (highest-rank objective method) \\

\midrule
\multicolumn{2}{l}{\bf Branch-and-bound algorithm} \\
$Q$ & Stack of visited unexplored nodes \\
$I$ & Incumbent, i.e.\ lexicographically best-found solution \\
$u,v,r$ & 
Nodes ($r$ is the root) in the branch-and-bound tree \\
$\mathcal{S}(u)$ & Feasible solutions below node $u$ in the branch-and-bound tree \\
$\ell$ & Node level, i.e.\ job index, in the branch-and-bound tree \\
$t_i$ & Partial completion time of machine $M_i$ in a branch-and-bound node \\
$\mathcal{R}$ & Subset of jobs scheduled below a branch-and-bound node \\
$L_i,\vec{L}$ & Component $L_i$ of vectorial lower bound $\vec{L} = (L_1,\ldots,L_m)$ \\
$U_i, \vec{U}$ & Component $U_i$ of vectorial upper bound $\vec{U} = (U_1,\ldots,U_m)$ \\ 
$\tau$ & Time point \\
$\tilde{p}_j$ & Piece of job $J_j$ \\
$\lambda, \Lambda$ & Amount of processing time load \\

\midrule
\multicolumn{2}{l}{\bf Makespan recovery problem} \\
$I_{\init}$ & Initial instance $(m_{\init},\mathcal{J}_{\init})$ \\
$I_{\new}$ & Perturbed instance $(m_{\new},\mathcal{J}_{\new})$ \\
$S_{\init}$ & Initial optimal schedule for $I_{\init}$ \\
$S_{\rec}$ & Recovered schedule for $I_{\new}$ \\
$S_{\new}$ & Optimal schedule for $I_{\new}$ \\
$\rho$ & Approximation ratio \\

\midrule
\multicolumn{2}{l}{\bf Uncertainty modeling} \\
$\mathcal{U}(f,k,\delta)$ & Uncertainty set \\
$f$ & Perturbation factor \\
$k$ & Number of unstable jobs \\
$\delta$ & Number of new machines \\
$C_{\max}^*(m,\mathcal{J})$ & Optimal objective value of makespan problem instance $(m,\mathcal{J})$ \\

\midrule
\multicolumn{2}{l}{\bf Binding recovery strategy} \\
$T$, $T_{\new}$ & Target makespan for instance $I$ and perturbed instance $I_{\new}$ \\
$T_{\new}$ & Target makespan for perturbed instance $I_{\new}$ \\
$\mathcal{M}'$ & Subset of machines \\
$m'$ & Number of machines in $\mathcal{M}'$ \\
$\mathcal{J}'$ & Subset of jobs \\
$\eta$ & Reduction of $p_j$ \\
$(\hat{m},\hat{\mathcal{J}})$ & Neighboring instance of $(m,\mathcal{J})$ \\
$\hat{p}_j$ & Processing time of job $J_j$ in $(\hat{m},\hat{\mathcal{J}})$ \\
$\mathcal{M}^s, \mathcal{M}^u$ & Stable machines $\mathcal{M}^s$ 
and unstable machines $\mathcal{M}^u = \mathcal{M}\setminus\mathcal{M}^s$ \\
$m^s,m^u$ & Number of stable ($m^s$) and unstable ($m^u$) machines \\
$\mathcal{J}_{\init}^s$ & Subset of stable jobs in $I_{\init}$ \\
$\mathcal{J}_{\new}^s$ & Subset of stable jobs in $I_{\new}$ \\
$F$ & Maximum processing time in $I_{\new}$ \\

\midrule
\multicolumn{2}{l}{\bf Flexible recovery strategy} \\
$\mathcal{J}^B,\mathcal{J}^F$ & Subset of binding ($\mathcal{J}^B$) and free ($\mathcal{J}^F = \mathcal{J}\setminus\mathcal{J}^B$) jobs \\
$\mathcal{J}_i^B$ & Binding jobs originally assigned to machine $M_i$ \\
$\mu_j$ & Machine executing job $J_j$ in $S_{\init}$ \\
$g$ & Limit on migrations of binding job \\

\midrule
\multicolumn{2}{l}{\bf Numerical results} \\
$\kappa$ & Phase transition parameter \\
$\kappa^*$ & Critical value of phase transition parameter  \\
$b$ & Number of bits for generating processing times \\
$q$ & Processing time seed \\
$\mathcal{U}$ & Discrete uniform distribution \\
$\mathcal{N}$ & Normal distribution \\
$W$ & Weighted value, i.e.\ weighted sum of objective functions \\
$Ub,Lb$ & Best-found incumbent ($Ub$) and lower bound ($Lb$) \\
$d_m, d_n$ & Number of machine ($d_m$) and job ($d_n$) disturbances \\
$W^N$ & Normalized weighted value \\
$W^*$ & Best computed weighted value \\
$C_{\max}^N$ & Normalized makespan \\
$C_{\max}^*$ & Best recovered makespan \\

\bottomrule
\end{longtable}
\end{center}

\end{appendices}
\end{document}